\documentstyle[12pt]{article}
\newcommand{\C}{1\!\!\!C}

\newcommand{\R}{I\!\!R}
\newcommand{\Sc}{{\cal S}}
\newcommand{\be}{\begin{eqnarray*}}
\newcommand{\ee}{\end{eqnarray*}}

      \date{}
\begin{document}

\setcounter{page}{1}


\begin{flushright}
TPH-207(rus) 1995

hep-th/9602002
\end{flushright}

\begin{center}
\large{{\bf Complex structure and solutions \\
of classical nonlinear equation \\
with the interaction $u^4_4$
}
\footnote{This paper is the first part of the
 project $\phi^4_4\cap M.$}
\medskip
\medskip
\\Edward  P. Osipov
\\            Department of Theoretical Physics
\\Institute for Mathematics
\\630090 Novosibirsk 90
\\RUSSIA}
\medskip
\\E-mail address (Relcom): osipov@math.nsk.su
      \date{}
\end{center}

      \date{}
\begin{abstract}
We consider the (real) nonlinear wave equation
$$\Box u + m^2 u + \lambda u^3  = 0, \quad m > 0, \quad
\lambda > 0, $$ on four-\-dimensional Minkowski space.
  Let $U(t),$  $W,$ and $S$ be the (nonlinear) operator of
dynamics and, respectively, the (nonlinear) wave and scattering
operators for this nonlinear wave equation. We introduce the
complex structure and show that the operators $U(t),$
$W,$ and $S$ define complex analytic maps on the space of
initial Cauchy data with finite energy. In other words, let
$R(\varphi, \pi) = \varphi + i\mu^{-1}\pi$ be the map of initial
data on the positive frequency part of the solution of the free
Klein-\-Gordon equation with these initial data. The operators
$RU(t)R^{-1},$  $RWR^{-1},$ and $RSR^{-1}$ are defined correctly
and are complex analytic on the complex Hilbert space
$H^1({\R}^3,\C).$ In particular, for $z(\alpha ) = \sum_{1\le n\le N}
\alpha_n z_{n},$ $z_{n} \in H^1({\R}^3,\C),$
$\alpha_n\in \C,$ $\langle RU(t)R^{-1} z(\alpha ), h \rangle,$
 $\langle RWR^{-1} z(\alpha ), h \rangle,$
and $\langle RSR^{-1} z(\alpha ), h \rangle$
are entire antiholomorphic functions
in $\alpha_n,$  $1\le n\le N. $

\end{abstract}
\newpage



\medskip
{\large 1. Introduction}
\medskip
\medskip

We consider the classic nonlinear wave equation
$$ \Box u+m^2 u+\lambda u^3=0,\;\; m>0,\;\;
\lambda>0 \eqno(1.1) $$
in fo\-ur-\-di\-men\-si\-on\-al
Minkowski spa\-ce-\-ti\-me.
We introduce the (global) complex structure
and prove
that
the wave and scattering (nonlinear) operators
 for this equation
 are complex analytic maps
 on the space
 of fi\-ni\-te en\-ergy Cauchy data.
This means that
the wave and scattering operators
define
 the complex analytic maps of  positive frequency parts
 given by initial data with finite energy.

This result allows us to consider
 the quantum field with the help
of its Wick kernel
 and to construct a
 bilinear form which is the solution of nonlinear quantum wave
equation with the cubic nonlinearity
in
fo\-ur-\-di\-men\-si\-on\-al spa\-ce-\-ti\-me
\cite{He74,Osi84,Osi94.1,OsiMIT94,Osi94.2}.
This bilinear form-\-so\-lu\-ti\-on
 is defined in the Fock space
of the free $in$-\-field on the subspace, generated by linear
 combinations of coherent vectors
 corre\-s\-pond\-ing to coherent
 $in$-\-states with finite energy.

To consider vacuum averages and the integral of the
type of
Paneitz, Pedersen, Segal, Zhou
\cite{PPSZ91}
we  need namely this holomorphity of Wick kernels.

The complex structure for the solutions of nonlinear
 equations is considered as follows.
We introduce the indexation of
 solutions of free equations with
 the help of complex variable.
This complex variable is the
 positive frequency part of a free real solution
 considered at a fixed time (for instance,
 at time zero).
 The unique (real) free
solution (with finite energy)
 can be defined uniquely by its
positive frequency part
 $\varphi^+,$ as well as by its initial data,
 (real) canonical coordinate and
 (real) canonical momentum.

In this
 case the map
 $$
\varphi^+:=\varphi+i\mu^{-1}\pi=R(\varphi,\pi),\quad
R^{-1}\varphi^+=(\mbox{ Re }\varphi^+,\mu\mbox{ Im
}\varphi^+),
$$
 where $(\varphi,\pi)\in H^1(\R^3)\oplus L_2(\R^3),$
$\varphi^+\in H^1(\R^3,\C),$
$\mu=(-\Delta+m^2)^{1/2},$
 defines the one-to-one
correspondence
between these variables.
Thus, we have a possibility to index solutions of nonlinear
 equation with the help
of complex variables.

The main consequence of our consideration is the following
principal result. The dynamics of the non-\-linear equation is
given by (infinite-\-dimensional) holomorphic maps. Namely, we
construct the complex structure for the nonlinear equation
(1.1) and we prove that the nonlinear dynamics, wave and
scattering operators generate
holomorphic maps of positive frequency parts
defined by solutions.


The proof of holomorphity of these maps is closely
connected with unitarity of corres\-pond\-ing operators.
 (these operators appear when we consider
corresponding quantum objects).
To proof holomorphity of these maps we
 consider the complex
structure for small initial data, i.e. for
 initial data from a
neighborhood
of zero in the space
 $H^1\oplus L_2.$ It means that we consider solutions with
small energy.
At this point, it is
very important
the stability of solutions. This stability of solutions gives
(for sufficiently small initial data) the uniqueness
of essential unitarizability (of the derivative) of
the wave operator and of the scattering operator.

The stability is closely related with the fact
that the equation corresponds physically to
a massive (scalar) field and
for small energies (less than
the double mass constant)
particles cannot be  created. In addition,
this nonlinear equation
has no solutions interpreted as bounded states.

Then, the consideration of  solutions smoothed with appropriate
functions (= functions with compact support in momentum space)
allows to extend the complex analyticity on large initial data.
Moreover, this complex analyticity due to the real analyticity
over (real) initial data gives the values that are equal to the
values  given by solutions with large initial data.

This  global complex structure gives  entire holomorphic
functions on finite-\-dimensional subspaces.
 With the help of uniform convergence we can easily
extend this complex
structure on larger class of smoothing
functions and on the all solutions
with initial data with finite
energy.
The (anti)holomorphity of functions
 $\langle RWR^{-1}(z_{in}(\alpha)),h \rangle$ on
finite-\-dimensional subspaces and integrals over these
subspaces
(see
Paneitz, Pedersen, Segal, Zhou
\cite{PPSZ91}
are very   important
for the construction of the quantum field.
 This holomorphity  with the help
 of Wick kernels
 allows
to construct the quantum field and its vacuum averages
(see
 Heifets
\cite{He74},
 Osipov
\cite{Osi94.1,OsiMIT94,Osi94.2}).

  Thus, we prove
the following theorem.

\medskip
{\large Theorem 1.1.}
  {\it Let $R(\varphi,\pi)=\varphi+i\mu^{-1}\pi$
 be the isomorphism
$(H^{1/2} \oplus H^{-1/2}, J)$
 onto
$H^{1/2}(\R^3,\C).$
 Let $W$
 be the wave operator for the nonlinear wave
equation $(1.1).$
 Then the map $RWR^{-1}$
is defined correctly as
a map of $D(\mu^{1/2})$
 on $D(\mu^{1/2})$
 and is a complex analytic map
of the complex Hilbert space
 $D(\mu^{1/2})$
$\simeq$ $ H^1(\R^3,\C)$
 onto itself.  In particular,
for $z_{in}(\alpha) = \sum^N_{j=1}\alpha_j
z_{in,j},$ $z_{in,j} \in D(\mu^{1/2}),$ $ \alpha_j \in \C$
 (i.e.
$z_{in,j} \in H^1(\R^3,\C)),$ $h \in H^{1/2}(\R^3,\C$),
 the functions $\langle RWR^{-1} z_{in}(\alpha),
h\rangle_{H^{1/2}(\R^3,\C)}$
 are entire
 antiholomorphic
functions on
$(\alpha_1,...\alpha_N) \in {\C}^N.$

 The same assertion is valid for
the map $RSR^{-1},$   where $S$
 is the scattering operator
of the nonlinear wave equation $(1.1).$}

\medskip
{\it Remark.}
  We emphasize that  variables
$(\varphi,\pi),$  $\varphi^+ = \varphi + i\mu^{-1}\pi$
should be considered as vectors from the spaces
 $(H^{1/2}\oplus H^{-1/2}, J),$ $H^{1/2}(\R^3,\C),$
respectively.
The complex structure corresponding to the
 symplectic 2-\-form of the nonlinear wave equation acts
in the complex Hilbert space
 $(H^{1/2}\oplus H^{-1/2}, J),$ namely.
 The vector $\varphi^+$ from $H^{1/2}(\R^3,\C)$
 defines  the  (free real) solution
$\exp(i\mu t)\varphi^++\exp(-i\mu t)\varphi^-.$
 In momentum space this vector $\varphi^+$
is defined by the
  function
$\mu(k)(\varphi^+)^\sim(k).$
Here
$(\cdot)^\sim$ is the Fourier transformation
and the   function
$\mu(k)(\varphi^+)^\sim(k)$  is
a square integrable function
 over the measure
$\vartheta(k_0)\delta(k^2-m^2)d^4 k$ (in momentum space).

At the same time the solution of the nonlinear equation have
the initial data from
$H^1\oplus L_2$
and the natural map
$$
R(\varphi,\pi)=\varphi+i\mu^{-1}\pi$$
 is the isomorphism
 $H^1\oplus L_2\to
H^1(\R^3,\C).$
 Note that the space $H^1(\R^3,\C)$ is the domain
 $D(\mu^{1/2})$ of definition
of the operator  $\mu^{1/2}$
in the space $H^{1/2}(\R^3,\C).$
This natural map
appears as the change of initial data
on  positive frequency
 part of the same solution
 (or as the change of a pair  of real
variables on the unique complex variable).

  Therefore, a unique solution is represented
 uniquely by the (unique) pair of initial data
$(\varphi,\pi),$ by the pair of real
 variables $(\varphi, \mu^{-1}\pi),$
 by the positive frequency part of the solution
$\varphi^+(t,x),$
 by the complex variable
 $$ \varphi^+(0,x) = \varphi(x) + i (\mu^{-1}\pi)(x),$$
 and by the
square integrable
 (over  the measure ${d^3 k / \mu(k)}$) complex
function $\mu\varphi^\sim(k) + i\pi^\sim(k).$ The  last
complex-\-valued function  is used in the isomorphism
$H^1\oplus L_2 \to L_2(\R^3,\C)$
defined by
Baez, Zhou
\cite{BaezZhou90,Baez91-92}.
This isomorphism gives the complex variable,
corresponding to the solution ( in this case
this variable is
 $\varphi(x) +
i\mu^{-1}\pi(x)$
 in spatial coordinates  and
$\mu(k)\varphi^\sim(k) + i\pi^\sim(k)$ in momentum
coordinates).

In addition, the transformation
$$
(P(\varphi,\pi))^\sim=\mu\varphi^\sim(k)+i\pi^\sim(k)
$$
is the projection on the massive hyperboloid by Goodman
 \cite{Goodman64}.
 That is, by definition
\be
 \int v^+(t,x)
f(t,x)dtdx &=& \int \mu^{1/2} v^{+}(0,x)
 \mu^{1/2}\mu^{-1} \int
e^{i\mu t}f(t,x)dt d^3x \cr =
\int v^{+\sim}(0,k) \int e^{i\mu t}
({\cal F}_3f) (t,k)dtd^3 k &=& \int v^{+\sim}(0,k)
f^\sim(\mu,k)d^3 k \\
= \int \mu v^{+\sim}(0,k) f^\sim(\mu,k)
{d^3k\over\mu} &=& \langle v,Pf \rangle_{H^{1/2}(\R^3,\C)},
\ee
$$(Pf)(x) = \mu^{-1} \int e^{i\mu
t}f(t,x)dt$$
and $$(Pf)^\sim(k)=
\int e^{i\mu t}e^{ikx}f(t,x)dtd^3 x.$$

Thus, if a solution $u$ has the initial data
$(\varphi,\pi),$ then its positive
frequency part is $\varphi^+,$ $\varphi^+(x)=
\varphi(x)+i(\mu^{-1}\pi)(x),$ and the projection
on the massive
hyperboloid in momentum space is equal to
 $$ (P(u))^\sim(k)
=e^{i\mu(k)t}(\mu(k)\varphi^\sim(k)+i\pi^\sim(k)). $$
In this case the solution itself is equal to
 \be u(t,x) &=&(\cos\mu t \varphi)(x)+({\sin\mu
t\over\mu} \pi)(x) \\ &=& (\exp i\mu t
\varphi^+)(x)+(\overline{\exp i\mu t \varphi^+})(x)
\\ &=&
(\cos\mu t\mbox{ Re}(\varphi^+))(x) +
(\sin\mu t\mbox{ Im}(\varphi^+))(x)
\\
&=&
\int e^{i\mu t} e^{ikx} P(u)^\sim(k) {d^3 k\over
2\mu(k)}
\\
&=&
\int e^{ik_0 t+ikx} P(u)^\sim(k_0,k)
\vartheta(k_0)\delta(k^2-m^2)d^4 k.
\ee

Therefore, the isomorphism introduced by
Baez, Zhou
 \cite{BaezZhou90,Baez91-92},
and Osipov
 \cite{Osi94.2}
agrees with the isomorphism introduced in our paper.
In this case an isomorphism maps the same solution
in the different spaces with the help of which we
index solutions: that is the space $H^1\oplus L_2$ for
initial data given by vector $(\varphi(x),\pi(x))$ in
coordinate space,
or $(\varphi^\sim(k),\pi^\sim(k))$ in
momentum space,
or the space $H^1({\R}^3,\C)$
 for positive
frequency part   $\varphi^+,$ $\varphi^+ =
\varphi + i\mu^{-1}\pi,$ of solutions with finite
energy, or
 the Baez-\-Zhou isomorphism
 \cite{BaezZhou90,Baez91-92},
or the projection by
Goodman
 \cite{Goodman64}
on the massive hyperboloid given by
 the vector
 $\mu\varphi^\sim(k)+i\pi^\sim(k)$
 in momentum space with
 the measure $\vartheta(k_0)\delta(k^2-m^2)d^4 k$
 (= a complex valued square integrable function
 in momentum space with the measure
 $d^3 k/2\mu(k)$).



\medskip
\medskip
\medskip
\medskip
{\large 2. Global complex structure.
 Ideas of the proof of Theorem 1.1
}
\medskip
\medskip

In this section we outline ideas and principles
of the proof of Theorem 1.1. The proof of
Theorem 1.1
 can be reduced to the proofs of theorems that are given
in the next sections. We complete Section 2
by the deduction of Theorem 1.1 from Theorem 5.1.

First of all, we
 note that the isomorphism
 $R$
maps the space of initial data $H^1 \oplus L_2$
of solutions with finite energy  into the space
 $H^1(\R^3,\C)$
of positive frequency parts of these solutions.
The complex Hilbert space $H^1(\R^3,\C)$
is identified with the domain $D({\mu}^{1/2})$
of definition of the operator
$\mu^{1/2}$  with topology
given by the norm
$$\Vert\mu^{1/2}(\cdot) \Vert_{H^{1/2}(\R^3,\C)}.$$
$R^{-1}$ is correctly defined on
$D(\mu^{1/2})$ and  $$R^{-1}D(\mu^{1/2})
= H^1 \oplus L_2\simeq
L^1_2\oplus L_2.$$

If $H_1,$ $H_2$ are complex Hilbert spaces and
$F$ is the map of $H_1$ into $H_2,$
then to prove the complex holomorphity of the map
$F$ it is sufficient to proof that for every
$x, h \in H_1,$ $y \in H^\ast_2$ ($H^\ast_2$
the (complex)
dual to $H_2,$ $H^\ast_2 \simeq H_2$)
any function $(y,F(x+\alpha h))_{H_2}$
is  holomorphic in $\alpha \in \C$ for
sufficiently small $|\alpha|,$
see
\cite[ch. 2.3, p.84-85]{Berger77}.

The basic idea of the proof of Theorem 1.1 is  the
following.

The first step  is the proof
of holomorphity of maps $RWR^{-1}$ and  $RSR^{-1}$
in a neigh\-bor\-hood of small initial data
with finite energy (see Section 3).
Since a real analyticity is proved by
Kum\-lin
\cite{Kumlin92},
see also
  Appendix and
\cite{BaezZhou89,Baez91-92,He74,RaczkaStr79,MorStr72,MorStr73},
so to prove
 the complex holomorphity we need to prove
commutativity of derivatives of wave operator and of
scattering operator with the operator of imaginary
unite.
To obtain
this result,
 we prove
for our purpose
the uniqueness of essential unitarizability for the
first derivatives of the wave operator and of the
scattering operator in some small zero neighborhood
of (good) initial data (Theorem 3.1).
For this purpose we
use Krein and coauthors ideas
\cite{DaletskiiKrein71}
and results of
 Pa\-neitz  and Se\-gal
 \cite{Pan80-82,Pan81,Pan82,PanSe80}.
 These  ideas and results give a possibility
to use a stability of nonlinear solutions and prove the
uniqueness of essential unitarizability
for the first derivative of wave
equation and for the first derivative of scattering operator.
As good initial data we use
 the following initial $in$--data
 $u_{in}.$ The initial $in$--data is such that
$u_{in}$ and $Ju_{in}$
(where  $J$ is the operator of im\-ag\-in\-ary unit,
$J=R^{-1}iR$) belong to the
intersection (of finite number) of some Hilbert
and Ba\-nach spaces
(and are dense in these spaces).
These, in general, various spaces of initial data
were considered and used by
Mo\-ra\-wetz, Stra\-uss
\cite{MorStr72,MorStr73},
 Se\-gal
\cite{Se66,Se68},
Pa\-neitz
\cite{Pan80-82}.
 The intersection of these spaces contains the
Schwartz space
${\Sc}_{Re}(\R^3)\oplus{\Sc}_{Re}(\R^3).$
The treatment of finite-\-dimensional zero
neighborhood of this space is sufficient for our
constructions.

This result about the unitar\-iz\-abil\-ity yields
us the commutativity of operator of im\-ag\-in\-ary
unit
with the derivative of the wave operator and with the
derivative of scattering operator
of the nonlinear equation depending on small initial
$in$--data.

This commutativity with the operator of
 im\-ag\-in\-ary unit allows us to obtain the complex
holomorphity on a zero neighborhood
of some fi\-ni\-te-\-di\-men\-si\-on\-al
subspace of initial
$in$-\-data. This
set of finite-\-dimensional subspace
is dense in  the mentioned above (finite) intersection
of Banach spaces of initial $in$--data. Instead  of the
mentioned above (finite) intersection of Banach spaces of
initial $in$--data we may choose and we choose
the Schwartz space
${\Sc}_{Re}(\R^3)\oplus{\Sc}_{Re}(\R^3).$
The $C^\infty$ Frech{\'{e}}t differentiability of nonlinear
solution and the complex holomorphity on
finite-\-dimensional
zero  neighborhoods allow us to extend the commutativity
of imaginary unit on all derivatives of wave operator and
scattering operator at zero.

This result and the real analyticity yield us the complex
 analyticity for initial
$in$--data in
   a zero neighborhood in $H^1(\R^3,\C),$
that is, for $u^+_{in}\in H^1(\R^3,\C),$
$\Vert u_{in}^+\Vert_{H^1(\R^3,\C)} < \vartheta$
for some $\vartheta>0$ (Theorem 3.2). Here
$u^+_{in}$ is the positive frequency part of free
solution (with corresponding initial data).

In order to extend this holomorphity on large initial data we
consider the functions  that are the smoothed (positive
frequency part of nonlinear) solution. These smoothed
solutions are smoothed over temporal
and spatial coordinates with appropriate (rapidly decreasing
in coordinate space) test functions with compact support in
momentum space.

Due to a (strict) positivity of mass constant,
 real analyticity
(in the whole space), and complex holomorphity in a zero
neighborhood the smoothed solution coincides with a smoothed
polynomial on initial $in$--data (Theorem 4.1).

The
degree of this polynomial is bounded by the size
of test function support (in
momentum space)  and is less than
$$
\sup\{ [p_0/\mbox{{\it mass constant}}] + 1\; | \;
(p_0,{\bf p})\in\mbox{ supp } f^\sim\}.
$$

Thus, we have the smoothed solutions and the smoothed
polynomials.
The smoothed solutions are real
analytic functions and the
smoothed polynomial are
holomorphic functions on the whole space of initial
data with finite energy. These considered smoothed
functions coincide on some zero neighborhood.
Thus, the uniqueness (Lemma  5.2), the
 real analy\-t\-i\-c\-ity
of  smoothed functions, the complex holomorphity of
smoothed polynomials give the coincidence on the
whole space
 $H^1(\R^3,\C).$ In other words, the smoothed solutions are
 complex analytic on the whole space of initial data with
finite energy.

Uniform continuity allows to extend the  holomorphity on
more wide class of test functions.
 This class of test  functions contains (smooth) functions
with compact support  (in  space-\-time  coordinates)
 too (Theorem 5.1). The equation of motion allows also to
extend the holomorphity
(of smoothed solutions) on test functions of the form
 $\delta(t)f(x), f\in\Sc(\R^3).$
This allows to complete the proof of Theorem 1.1.

Using outlined ideas we give  the detailed
 proof of
 main statements in Sections 3-5  of the paper. These main
statements are Theorems 3.1, 3.2, 4.1, 5.1.

In addition, in Sections 6, 7 we give more detailed
consideration of the assertions used
by
Paneitz and Segal
\cite{PanSe80,Pan80-82}
 for the proof of a uniqueness of essential
 unitarizability of operators $dW(u_{in})$
and $dS(u_{in})$ for small initial data.
These assertions are Theorems 6.1-6.3,
Corollary 6.4, and Theorem 7.1.

In Appendix we give also the improvement of
 Kumlin proof
\cite{Kumlin92}
of real analyticity for solutions of the
 nonlinear equation.

To complete this Section we deduce
Theorem 1.1 as a consequence
of Theorem 5.1.

\medskip
{\it Proof of Theorem $1.1.$}
  We derive
 Theorem 1.1 from Theorem 5.1.
Theorem 5.1 implies that the function
$$ \int\langle
RWU_0(t)R^{-1}(u^+_{in}(\alpha)),
h\rangle_{H^{1/2}(\R^3,\C)}
\chi(t)dt,
$$
$h\in\Sc(\R^3,\C),$ $\chi\in L_1(\R),$
$u^+_{in}(\alpha) =
\sum^N_{j=1} \alpha_j u^+_{in,j},$
 $u^+_{in,j} \in H^1(\R^3,\C),$
is an entire antiholomorphic function
in  $\alpha\in\C^N.$

Let
$\chi\in\Sc_{Re},$
$\mbox{ supp }\chi\subset[-1,1],$
$\chi(t) \ge 0,$ $\int\chi(t)dt = 1,$
$\chi_\sigma(t) = \sigma\chi(\sigma t).$
For   $\sigma\to\infty$
$\chi_\sigma(t)$ converges to the $\delta$--function
(in the sense of generalized functions).
It is clear that for
$\sigma\to\infty$
$$
\int\langle
RWU_0(t)R^{-1}(u^+_{in}(\alpha)),
h\rangle_{H^{1/2}(\R^3,\C)}
\chi_\sigma(t)dt
\quad
(\; = \int\langle u^+(t),h\rangle_{H^{1/2}(\R^3,\C)}
\chi_\sigma(t)dt) $$
converges to
 $\langle RWR^{-1}(u_{in}^+(\alpha)),
h\rangle (=\langle u^+(0),h
 \rangle\equiv\langle u^+,h\rangle).$

To obtain the uniform bound,
we use the equality
 $f(1)-f(0)=\int^1_0 ds{d\over ds}f(s),$ and write
$$
\int
dt\chi_\sigma(t)\langle
u^+(t),h\rangle_{H^{1/2}(\R^3,\C)}-
\langle u^+,h\rangle_{H^{1/2}(\R^3,\C)} $$
$$= \int
dt\chi_\sigma(t)(\langle u^+(t),
h\rangle_{H^{1/2}(\R^3,\C)}-
\langle u^+,h\rangle_{H^{1/2}(\R^3,\C)})
$$
$$=\int dt\chi_\sigma(t)t \;\int^1_0
ds\langle\dot u^+(st),h
\rangle_{H^{1/2}(\R^3,\C)}.
$$
The equation of motion
(in the sense of generalized functions
in spatial coordinates)
implies that
\be u^+(t) &=& u(t)+i\mu^{-1}\dot u(t) \cr \dot
u^+(t) &=&
\dot u(t)+i\mu u(t)+i\mu^{-1}\lambda u^3(t).
\ee
The Sobolev inequality
 and the energy conservation
imply  that
 \be
\Vert\dot u^+(t)\Vert_{L_2(\R^3,\C)}
&\leq & c_1(H(u(t))^{1/2} +
H(u(t))^{3/2}) \cr
&=& c_2(H_0(u_{in}(\alpha))^{1/2}+
H_0(u_{in}(\alpha))^{3/2})
\cr
&\leq& c(u_{in,1},..., u_{in,N})
(|\alpha|+|\alpha|^{3/2}).
\ee
Here
$H(u)$ is
the total energy of  solution $u$
and $H_0(u_{in})$
is the free energy of the (free) solution
$u_{in}.$
Thus, we obtain the estimate
$$ |\int dt\chi_\sigma(t)t\;
\int ds\langle\dot
u^+(st),h\rangle_{H^{1/2}
(\R^3,\C)}|\leq\sigma^{-1}\sup_{|t|\leq
1/\sigma}|\langle\dot u^+(t),h
\rangle_{H^{1/2}(\R^3,\C)}| $$ $$
\leq\sigma^{-1}c(u_{in,1},...,u_{in,N})
(|\alpha|+|\alpha|^{3/2})
\Vert h\Vert_{H^1(\R^3,\C)}
$$
and the uniform convergence
$$
\int dt\chi(t)\langle
RWU_0(t)R^{-1} u^+_{in}(\alpha),h\rangle_{H^{1/2}
(\R^3,\C)} \to\langle RWR^{-1}
u_{in}^+(\alpha),h\rangle_{H^{1/2}(\R^3,\C)}.  $$
Therefore, the function
$\langle RWR^{-1}
u_{in}^+(\alpha),h\rangle_{H^{1/2}(\R^3,\C)},$
$h\in\Sc(\R^3, \C),$
is antiholomorphic in $\alpha.$
If now
$h\in
H^{1/2}(\R^3,\C), h_n\in\Sc(\R^3,\C)$
and $h_n\to h$ in
$H^{1/2}(\R^3,\C),$
then $$\langle RWR^{-1} (u_{in}^+(\alpha)),
h_n\rangle_{H^{1/2}(\R^3,\C)}$$
converges uniformly to $$\langle
RWR^{-1} u_{in}^+(\alpha),
h\rangle_{H^{1/2}(\R^3,\C)}.$$
This is a consequence of the
 uniform estimate \be && |\langle
RWR^{-1}u_{in}^+(\alpha),
h_n\rangle_{H^{1/2}(\R^3,\C)} - \langle
RWR^{-1} u_{in}^+(\alpha),
h\rangle_{H^{1/2}(\R^3,\C)}| \cr &&
\leq
 \Vert
RWR^{-1}u^+_{in}(\alpha)\Vert_{H^{1/2}(\R^3,\C)}
\,\Vert h_n-h\Vert_{H^{1/2}(\R^3,\C)}
\cr
&&
\le c \Vert
u^+_{in}(\alpha)\Vert_{H^1(\R^3,\C)} \,\Vert
h_n-h\Vert_{H^{1/2}(\R^3,\C)}
\cr
&&\leq c(u_{in,1},...,
 u_{in,N})\,|\alpha|\,\Vert h_n-h\Vert_{H^{1/2}(\R^3,\C)}.
 \ee
Moreover, the following uniform estimate is valid
 $$ |\langle RWR^{-1}
u_{in}^+(\alpha),h\rangle_{H^{1/2}(\R^3,\C)}|
$$
$$
\leq
\Vert\mu^{1/2}
RWR^{-1} u_{in}^+(\alpha)\Vert_{H^{1/2}(\R^3,\C)}
\,\Vert\mu^{-1/2}h\Vert_{H^{1/2}(\R^3,\C)}
\leq c(u_{in})\,|\alpha|\,\Vert\mu^{-1/2}h\Vert_{H^{1/2}}.
$$
This implies that the function
$RWR^{-1} u^+_{in}(\alpha)$
is holomorphic as a map of
 $H^1(\R^3,\C)$ into $D(\mu^{1/2}).$
In particular, the function
 $\mu^{1/2} RWR^{-1} u_{in}^+ (\alpha)$
is holomorphic.

This means also that if we write
the smoothness as an integral
(with the dualization in
$L_2(\R^3,\C),$ or in $L_2(\R^4,\C)$),
then
$$ \int
u^+(t,x)h(x)dx=\int RWR^{-1}(u^+_{in}(\alpha))(x) h(x) dx $$
is defined correctly and is a holomorphic function
in $\alpha$ for all
$h \in H^{-1}(\R^3,\C),$ that is, in particular, for
 $h = \delta(x_1) h(x_2,x_3) $ also.

Hence, these functions are complex analytic
on finite-\-dimensional subspaces and the map
 $RWR^{-1}$ is complex holomorphic as a map from
 $D(\mu^{1/2})$ on $D(\mu^{1/2})$  in
$H^{1/2}(\R^3,\C).$

The scattering operator $S$ can be considered
analogously.
  {\it Theorem $1.1$ is proved.}

\medskip
{\it Remark.}
  The essential point
(the fine tuning of considered spaces!)
is that we need to consider
the holomorphity in the space
$H^{1/2}(\R^3,\C)$ namely.
This follows from symplecticity
 (and K{\"a}hlerian structure)
of the space with respect to the given symplectic 2-form.
The symplecticity, the existence (and the uniqueness)
 of complex structure and of  operator of
imaginary unit (with requirements of
Poincare-\-invariance)
have  a consequence, that the
 symplecticity implies the equalities
 $$ U^+JU=J $$
and if $(U^+)^{-1} = U$ (unitarity), then
$$ JU=UJ.  $$




\medskip
\medskip
\medskip
\medskip
{\large 3. Complex structure and
 complex  analyticity
in a zero neighborhood
}
\medskip
\medskip

In this Section we consider the complex structure
and complex analyticity for solutions of the nonlinear
equation with small energy,
i. e. for the initial Cauchy $in$--data from
 a zero neighborhood
in $H^1\oplus L_2.$

We introduce some notation that we need.
Let $W,$ $S$ be the nonlinear wave operator and the
nonlinear scattering operator, respectively, for the nonlinear
wave equation 1.1.
The nonlinear operators $W,$ $S$ are  $C^\infty$
Frech{\'e}t differentiable invertible
maps on the  (real) Hilbert space
$H^1\oplus L_2$ and these maps are
$C^\infty$ symplectomorphisms,
see, for instance
\cite{Baez91-92}
 and references  therein.

Let $\mu$ be the operator, $\mu:=(-\Delta+m^2)^{1/2},$ where
$\Delta$ is the Laplacian and $m,$ $m>0,$ is a mass
constant from Eq. (1.1).

Let $R$ be the operator of isomorphism.
This operator maps a pair of
real-\-valued functions (on $\R^3$) onto the
complex-\-valued  function,
$$ R(\varphi,\pi)=\varphi+i\mu^{-1}\pi := \varphi^+,
\quad
R^{-1}\varphi^+=(\mbox{ Re }\varphi^+,\mu\mbox{ Im }\varphi^+).
$$
Let $J = R^{-1}iR.$ Then $J$ is the
(orthogonal) operator of imaginary unit,
$$J (\varphi,\pi)
= (-\mu^{-1}\pi, \mu\varphi),\quad
J^2 = -I.$$
In particular,  $R$ defines
a real linear isomorphism
of the space
$H^{1/2} \oplus H^{-1/2}$ on
$H^{1/2}(\R^3,\C)$
and the complex linear isomorphism
$(H^{1/2} \oplus H^{-1/2}, J)$
on  $H^{1/2}(\R^3,\C).$
The inner product in
$(H^{1/2} \oplus H^{-1/2}, J)$ is given by the expression
 \be &&
\langle(\varphi_1,\pi_1),(\varphi_2,\pi_2)
\rangle_{(H^{1/2}\oplus H^{-1/2},J)} \cr &&
=((\varphi_1,\pi_1),(\varphi_2,\pi_2))_{H^{1/2}\oplus
H^{-1/2}}
+i(J(\varphi_1,\pi_1),
(\varphi_2,\pi_2))_{H^{1/2}\oplus H^{-1/2}}
\cr
&&
=\omega((\varphi_1,\pi_1),J(\varphi_2,\pi_2))
+i\omega((\varphi_1,\pi_1),(\varphi_2,\pi_2)),
\ee
where $\omega,$ the imaginary part of the inner product,
 is the symplectic 2-form,
$$
\omega((\varphi_1,\pi_1),(\varphi_2,\pi_2))=
\int d^3 x (\varphi_1(x)\pi_2(x)-\pi_1(x)\varphi_2(x)).
$$
This choice corresponds to the choice of
the inner product in the form
 $$
 \langle(\varphi_1,\pi_1),(\varphi_2,\pi_2)
\rangle_{(H^{1/2}\oplus H^{-1/2},J)} =
\langle\varphi^+_1,\varphi^+_2\rangle_{H^{1/2}(\R^3,\C)}
= \int d^3 x
\overline{(\mu^{1/2}\varphi^+_1)}(x)
(\mu^{1/2}\varphi^+_2)(x), $$
i.e. the inner product is antilinear in the first
argument.

By $d^n F$ we denote Frech{\'e}t derivatives
of the transformation $F$ of a Banach space into an other
Banach space.

\medskip
{\large Theorem 3.1.}
  {\it Let  $u_{in}(\alpha)=\sum^N_{j=1} \alpha_j u_{in,j},
$ $\alpha = (\alpha_1,...,\alpha_N)\in\C^N,
$ $u_{in,j}\in\Sc_{Re}(\R^3)\oplus\Sc_{Re}(\R^3).$
Let $u(\alpha)$
be the solution of the nonlinear wave equation
$(1.1)$ with the initial $in$--data
$u_{in}(\alpha).$

There exists  a $\vartheta>0$ (depending on
$u_{in,j}, j=1,...,N$), such that for $|\alpha|<\vartheta$
the operators $dW(u_{in}(\alpha)), dS(u_{in}(\alpha))$
are correctly defined on $H^1\oplus L_2$
and are   bounded operators (on $H^1\oplus L_2$). These
bounded operators are uniquely essentially unitarizable
in the complex Hilbert space
$(H^{1/2}\oplus H^{-1/2},J),$ and
on $D(\mu^{1/2})$
$(=(H^1\oplus L_2,J)\subset (H^{1/2}\oplus H^{-1/2},J)) $
$$
\langle dW(u_{in}(\alpha))v_{in,1},\;
dW(u_{in}(\alpha))v_{in,2}
\rangle_{(H^{1/2}\oplus H^{-1/2},J)}
=\langle v_{in,1},v_{in,2}
\rangle_{(H^{1/2}\oplus H^{-1/2},J)},
$$
$$
JdW(u_{in}(\alpha))v_{in}=dW(u_{in}(\alpha))Jv_{in},
$$
$$
\langle dS(u_{in}(\alpha))v_{in,1},\;
dS(u_{in}(\alpha))v_{in,2}
\rangle_{(H^{1/2}\oplus H^{-1/2},J)}
=\langle v_{in,1},v_{in,2}
\rangle_{(H^{1/2}\oplus H^{-1/2},J)},
$$
$$
JdS(u_{in}(\alpha))v_{in}=dS(u_{in}(\alpha))Jv_{in}.
$$
By  continuity these equalities for linear operators
can be  extended uniquely
(as linear operator on the variable
$v_{in}$) on the whole complex space
 $(H^{1/2}\oplus H^{-1/2},J).$
}

\medskip
{\it Remarks.}
  1. The assertion of Theorem 3.1 on the argument
 $u_{in}$ can be extended by continuity on the zero
neighborhood of some Banach space of initial data,
or on {\it finite} intersection of Banach space
used by I. Segal
\cite{Se66,Se68},
 Morawetz and Strauss
\cite{MorStr72,MorStr73},
Heifets
\cite{He74},
Raczka and Strauss
\cite{RaczkaStr79},
 Paneitz
\cite{Pan80-82,Pan81,Pan82},
Kumlin
\cite{Kumlin92}
 for  initial data and initial
$in$-\-data.
In particular,
Paneitz
\cite{Pan80-82},
 Paneitz and Segal
\cite{PanSe80}
 used
some Banach space of initial data
to prove an analogue of Theorem 3.1.
It is obvious that the Schwartz space
$\Sc_{Re}(\R^3)\oplus\Sc_{Re}(\R^3)$ is contained in and
is dense in these Banach spaces, has
a  stronger topology, and is a nuclear space.

2. We remark that in Theorem 3.1 (and in the next ones)
 $\vartheta$ denotes a small strictly positive
constant.
This constant $\vartheta$ gives a choice of (some)
neighborhood of the zero solution,
this choice depends on the topology that we need take into
account in the considered theorem.  In general, these
constants
are various.

3. It is interesting to consider
complex structure in many-\-dimensional space-\-time.
The stability for small initial data (with a non-\-zero mass)
is similar
\cite{PanSe80},
and at the same time the consideration of complex structure
opens new possibilities for
global solutions with large initial data.
But there is reason to think  that the usual
Schwartz space is not adequate to describe global
solution with large initial data.

\medskip
{\it Proof of Theorem 3.1.}
  It is clear that the Frech{\'e}t derivative of operators
 $W,$ $S$ is defined correctly  as a linear operator
on the space $H^1\oplus L_2$ and can be  extended by
continuity as a linear operator
(depending on the point $u_{in}(\alpha)$)
on more wide (Banach)
space.

To prove Theorem 3.1 we apply Corollary 2.3
\cite{PanSe80},
 Theorem 7
\cite{Pan81},
 Theorem 16.3 \cite[Theorem 16.3]{Pan82},
Corollary
\cite[Corollary, p.115]{Pan80-82},
 Theorem 6.3
\cite{Pan82}
and verify the required conditions.

In the case of the operator $S$
it is convenient to use directly
\cite[Corollary, p.115]{Pan80-82},
 see also
\cite[Theorem 4]{PanSe80},
 Theorem 1, 2
\cite[pp.114-115]{Pan80-82},
 Theorem 6.3
\cite{Pan82}.
In this case
Corollary
\cite[Corollary, p.115]{Pan80-82}
 implies directly the assertion of Theorem 3.1
for the initial
$in$-data from $\Sc_{Re}(\R^3)\oplus\Sc_{Re}(\R^3)$.
In this case  the conditions
 $$
(m^2-\Delta)^{5/4}u_{in}(t,\cdot),\;\;
(m^2-\Delta)^{3/4}\dot u_{in}
(t,\cdot)\in L_1(\R^3)
\quad\mbox{for some}\quad t$$
are fulfilled.

The condition
$$
\int^{+\infty}_{-\infty}
\Vert u(t,x)\Vert^2_\infty dt<2m
$$
is fulfilled for the choice
$u_{in}\in{\cal F},
\Vert u_{in}\Vert_{{\cal F}} <\vartheta$
for sufficiently small positive
 $\vartheta,$ see  Theorem 1(b)
\cite{MorStr72},
here ${\cal F}$ is the space of initial data, defined
by Morawetz and Strauss
\cite{MorStr72}.

Thus, there exists a constant $\vartheta,$ such that
$u_{in}(\alpha), |\alpha|<\vartheta,$ satisfy
all conditions of Corollary 2.3
\cite{PanSe80},
Theorem 7
\cite{Pan81},
Theorem 16.3
\cite[Theorem 16.3]{Pan82},
Corollary
\cite[Corollary, p.115]{Pan80-82}
and Theorem 4.4B
\cite[p. 110]{Pan80-82}
and
\cite[Corollary 4.4B, p. 491]{Se68}.

In the case of the wave operator $dW(u)$ we apply
Theorem 7
\cite{Pan81},
 Theorem 6.3
\cite[Theorem 6.3]{Pan82},
 and Theorems 1, 2
\cite{Pan80-82}
 for the interval $(-\infty,0]$
(or for the interval $[0,\infty)$ in the case
of the wave operator into
the future, we may also use intervals
$(-\infty,T]$ and $[T,\infty)$).

 The proof of absence of
``bounded states"  was given by Paneitz
\cite[Cor\-ol\-l\-ary, p. 115-116]{Pan80-82},
 see also
\cite{PanSe80}.
 We give  this proof in Section 7.
  {\it Theorem $3.1$ is proved.}

\medskip
We now  go  to the consideration of holomorphity.

\medskip
{\large Theorem 3.2.}
 {\it There exists a strict positive $\vartheta$
such that for $z\in D(\mu^{1/2})$
( in $H^{1/2}(\R^3,\C)$),
$$\Vert z\Vert_{H^1(\R^3,\C)} < \vartheta
$$ i.e. for
$$\Vert\mu^{1/2}z\Vert_{H^{1/2}(\R^3,\C)} < \vartheta ,
$$
the transformations
$RWR^{-1}(z),$ $RSR^{-1}(z)$ are holomorphic in
$z$ as transformations of open set  $\Vert
z\Vert_{H^1(\R^3,\C)}<\vartheta$ of the space $H^1(\R^3,\C)$
into the space
 $H^1(\R^3,\C).$ In particular, $$
RWR^{-1}(z)=\sum{1\over n!}d^n RWR^{-1}(0) (z,...,z),
$$
$$
RSR^{-1}(z)=\sum{1\over n!}d^n RSR^{-1}(0) (z,...,z),
$$
 the series converge uniformly for  $\Vert
z\Vert_{H^1(\R^3,\C)}<\vartheta,$ and
the functions $$ \langle
RWR^{-1}(z),h\rangle_{H^{1/2}(\R^3,\C)}, \quad\langle
RSR^{-1}(z),h\rangle_{H^{1/2}(\R^3,\C)}, $$ $h \in
H^{1/2}(\R^3,\C),$ are antiholomorphic in $z$
 for $\Vert z\Vert_{H^1(\R^3,\C)} < \vartheta. $
}

\medskip
\medskip
{\it Remark.}
  The value of  $\vartheta$ is defined by the radius of
convergence of the Taylor series at zero in the real Hilbert
space $H^1\oplus L_2$ and its topology.

\medskip
{\it Proof of Theorem $3.2.$}
  First of all, we note that it is possible to write the equality
 $$
\langle RWR^{-1} u^+_{in},h\rangle_{H^{1/2}(\R^3,\C)}=
\langle Wu_{in},R^{-1}h\rangle_{{\cal H}_S}.
$$
Here we denote ${\cal H}_S: = (H^{1/2} \oplus
H^{-1/2},J),$ $u^+_{in}$ is the positive frequency part,
corresponding to the initial data of $u_{in}.$

$W(u_{in})$ is
$C^\infty$ Frech{\'e}t differentiable
function in $H^1\oplus L_2.$
Therefore the equality $$
({\partial\over\partial\alpha^{(0)}}
-i{\partial\over\partial\alpha^{(1)}})
\langle W(u_{in}+\alpha^{(0)}v_{in}+
J\alpha^{(1)}v_{in}), R^{-1}h
\rangle_{{\cal H}_S}
$$
$$
=\langle dW(u_{in}+\alpha^{(0)}v_{in}+
J\alpha^{(1)}v_{in})
v_{in}, R^{-1}h\rangle_{{\cal H}_S}
$$
$$
-i\langle dW(u_{in}+\alpha^{(0)}v_{in}+J\alpha^{(1)}v_{in})
Jv_{in}, R^{-1}h
\rangle_{{\cal H}_S}
$$
$$
=\langle dW(u_{in}+
\alpha^{(0)}v_{in}+\alpha^{(1)}Jv_{in})v_{in},R^{-1}h
\rangle_{{\cal H}_S}
$$
$$
+ \langle J dW(u_{in}+
\alpha^{(0)}v_{in}+\alpha^{(1)}Jv_{in})Jv_{in},R^{-1}h
\rangle_{{\cal H}_S}
$$
is defined correctly for
$u_{in}, v_{in} \in H^1\oplus L_2,$
$\alpha^{(0)}, \alpha^{(1)} \in \R.$
The transformation $W$  is
 real analytic in the Hilbert space
 $H^1\oplus L_2$ and
 the transformation $RWR^{-1}$ is
real analytic in the Hilbert space
 $H^1(\R^3,\C)$
 (see
\cite{Kumlin92},
and also
\cite{He74,RaczkaStr79,BaezZhou89,Osi94.1,Osi94.2}).
 Let $u_{in} (\alpha) =
 \sum^N_{j=1} (\alpha^{(0)}_j +
\alpha^{(1)}_jJ) v_{in,j}$
( $= R^{-1} u^+_{in}(\alpha)),$ $v_{in,j} \in
{\Sc}_{Re}(\R^3)\oplus {\Sc}_{Re}(\R^3),$
$\alpha^{(0)}_j,
\alpha^{(1)}_j\in\R.$ Theorem 3.1
implies the existence
of  $\vartheta(u_{in},N)>0$
such that for $|\alpha|<\vartheta (u_{in},N)$
\be
\langle dW(u_{in}(\alpha))Jv_{in},
R^{-1}h\rangle_{{\cal H}_S}
&=&
 \langle J dW(u_{in}(\alpha))v_{in},
R^{-1}h\rangle_{{\cal H}_S}
\cr
&=& -i\langle dW(u_{in}(\alpha))v_{in},
R^{-1}h\rangle_{{\cal H}_S}.
\ee
Therefore,
$$
({\partial\over\partial\alpha^{(0)}}
-i{\partial\over\partial\alpha^{(1)}})
\langle W(u_{in}(\alpha)),
R^{-1}h\rangle_{{\cal H}_S}=0
$$
and the function
$\langle W(u_{in}(\alpha)), R^{-1}h\rangle$
is antiholomorphic. The holomorphity and
Cauchy-\-Riemann conditions
imply that $(k_j=0,1)$ \be &&\langle d^n
W(u_{in}(\alpha))(J^{k_1}v_{in,1},...,J^{k_n}v_{in,n}),
R^{-1}h\rangle_{{\cal H}_S}
\cr
&&=
(\prod^n_{j=1}{\partial^{k_j}\over\partial\alpha^{(k_j)}_j})
\langle W(u_{in}(\alpha)), R^{-1}h\rangle_{{\cal H}_S}\cr
&&=\langle J^{k_1+...+k_n}d^n W(u_{in}(\alpha))
(v_{in,1},...,v_{in,n}),R^{-1}h\rangle_{{\cal H}_S}
\cr
&&=
(\prod^n_{j=1}(-i)^{k_j}
{\partial^{k_j}\over\partial\alpha^{(k_j)}_j})
\langle W(u_{in}(\alpha)), R^{-1}h\rangle_{{\cal H}_S}.
\ee
Since
$W \in C^\infty (H^1 \oplus L_2,H^1 \oplus L_2),$
 so
 $d^n W(0)$ is a $n$-linear symmetric continuous bounded
operator. The density of ${\Sc}_{Re} \oplus {\Sc}_{Re}$
in $H^1\oplus L_2$ allows to extend the equality $$ d^n
W(0)(J^{k_1}v_{in,1},...,J^{k_n}v_{in,n}) =
J^{k_1+...+k_n}d^n W (0)(v_{in,1},...,v_{in,n}) $$ by
continuity in  $(v_{in,1},...,v_{in,n})$ on the whole space
$H^1\oplus L_2.$

Now, if  $\Vert u_{in}\Vert_
{H^1\oplus L_2}<\vartheta,$ where
$\vartheta>0$ and less than the radius of
 convergence of Taylor series at zero for the real expansion
in  $H^1\oplus L_2,$
see
\cite{Kumlin92},
(and also
\cite{He74,RaczkaStr79,BaezZhou89,Baez91-92,Osi94.1,Osi94.2}),
 for
$v_{in}\in H^1\oplus L_2$
and for sufficiently small
$\alpha^{(0)},\alpha^{(1)}\in\R,$
 i.e. for $$ \Vert u_{in} + \alpha^{(0)}v_{in} +
\alpha^{(1)}Jv_{in} \Vert_{H^1\oplus L_2}
< \vartheta, $$ we have
$$
\langle W (u_{in} +\alpha^{(0)}v_{in}+
\alpha^{(1)}J v_{in},R^{-1} h
\rangle_{{\cal H}_S}
$$
 $$=\sum^\infty_{n=1}{1\over n!}
\langle d^n W(0)
(u_{in}+\alpha^{(0)}v_{in}+\alpha^{(1)}J v_{in}, ...,
u_{in}+\alpha^{(0)}v_{in}+
\alpha^{(1)}J v_{in}), R^{-1} h
\rangle_{{\cal H}_S}
$$
$$
=\sum_{n_1+n_2+n_3=n}{1\over n_1! n_2! n_3!}
\langle d^n W(0)
(\underbrace{u_{in},...,u_{in}}_{n_1},
\underbrace{\alpha^{(0)} v_{in},...,
\alpha^{(0)} v_{in}}_{n_2},
$$
$$
\underbrace{\alpha^{(1)}Jv_{in},...,
\alpha^{(1)}Jv_{in}}_{n_3}),
R^{-1}h\rangle_{{\cal H}_S}
$$
$$
=\sum_{n_1+n_2+n_3=n}{1\over n_1! n_2! n_3!}
\langle\alpha^{(0)n_2}
(\alpha^{(1)}J)^{n_3}d^n W(0)
(\underbrace{u_{in},...,u_{in}}_{n_1},
\underbrace{v_{in},...,v_{in}}_{n_2+n_3},
R^{-1}h\rangle_{{\cal H}_S}
$$
$$
=\sum_{n_1+n_2=n}{1\over n_1! n_2!}
(\alpha^{(0)}-i\alpha^{(1)})^{n_2}\langle
d^n W(0)(\underbrace{u_{in},...,u_{in}}_{n_1},
\underbrace{v_{in},...,v_{in}}_{n_2},
R^{-1}h\rangle_{{\cal H}_S}.
$$
Since the series converges
absolutely for $\Vert u_{in}\Vert_{H^1\oplus
L_2} <\vartheta$ and $\Vert u_{in} + (\alpha^{(0)}
+\alpha^{(1)}J) v_{in}\Vert <\vartheta,$
i.e. for $|\alpha|
<\Vert v_{in}\Vert^{-1}
(\vartheta -\Vert u_{in}\Vert),$
so its sum is antiholomorphic in
 $\alpha^{(0)} + i\alpha^{(1)}.$
Therefore, the function
 $\langle RWR^{-1} (u^+_{in} +\alpha v^+_{in}),
h\rangle_ {H^{1/2}(\R^3,\C)}$
is antiholomorphic in $\alpha$
for sufficiently small $|\alpha|.$

In order to show that  $RWR^{-1}$ is holomorphic as a
transformation in
 $H^1(\R^3, \C),$ it is sufficient to show
(see
\cite[pp. 84-85]{Berger77}),
that for
 $h\in H^1(\R^3, \C),$
$\Vert u^+_{in}\Vert_{H^1(\R^3, \C)}
<\vartheta$ (i.e. $\Vert\mu^{1/2}
u^+_{in}\Vert_{H^{1/2}(\R^3, \C)}< \vartheta$),
the function
$$\langle RWR^{-1} (u^+_{in} +\alpha v^+_{in}),
h\rangle_{H^{1/2}(\R^3,\C)}$$
is an\-ti\-ho\-lo\-mor\-phic in $\alpha.$
It is evident that for $h\in \Sc(\R^3, \C)$
the expressions $$
\langle RWR^{-1} (u^+_{in} +\alpha v^+_{in}), h\rangle_{H^1
(\R^3,\C)}
=
\langle\mu^{1/2} RWR^{-1} (u^+_{in} +\alpha v^+_{in}),
\mu^{1/2}h
\rangle_{H^{1/2}(\R^3,\C)}
$$ $$
= \langle RWR^{-1} (u^+_{in} +\alpha v^+_{in}), \mu h
\rangle_{H^{1/2}(\R^3,\C)}
\eqno(3.1)
$$
 are correctly defined and an\-ti\-ho\-lo\-mor\-phic
in $\alpha$
(the function $\mu(\cdot)$ is real-\-valued and
the operator
$\mu$ commutes with the imaginary unit).
If now $h\in H^1(\R^3, \C)$
and
$h_n\in \Sc(\R^3, \C),$
$h_n\to h$ in $H^1(\R^3, \C),$ then the
uniform in $\alpha$
estimate implies that
\be &&\vert\langle RWR^{-1} (u^+_{in} +\alpha
v^+_{in}), h\rangle_{H^1 (\R^3,\C)}- \langle RWR^{-1}
(u^+_{in}
+\alpha v^+_{in}),
h_n\rangle_{H^1 (\R^3,\C)}\vert \cr &&\le
\Vert\mu^{1/2}RWR^{-1}(u^+_{in}+
\alpha v^+_{in})\Vert_{H^{1/2}(\R^3,\C)}
\,\Vert\mu^{1/2}(h-h_n)\Vert_{H^{1/2}(\R^3,\C)}\cr
&&\le
\Vert \mu^{1/2}(u^+_{in}+\alpha
v^+_{in})\Vert_{H^{1/2}(\R^3,\C)}
\,\Vert\mu^{1/2}(h-h_n)
\Vert_{H^{1/2}(\R^3,\C)}\cr
&&\leq (c_1+c_2|\alpha|)\,
\Vert h-h_n\Vert_{H^1(\R^3,\C)}.
\ee
This estimate follows from the energy conservation
and Equality (3.1).
The estimate is  uniform in $\alpha$ and implies
the uniform convergence of the an\-ti\-ho\-lo\-mor\-phic
function $$
\langle RWR^{-1} (u^+_{in} +\alpha v^+_{in}),
h_n\rangle_{H^1(\R^3,\C)} $$ to
the function $$ \langle RWR^{-1}
(u^+_{in} +\alpha v^+_{in}), h\rangle_{H^1(\R^3,\C)} $$
and its antiholomorphity. Thus, the theorem is proved
for the case of the wave operator.

The operator  $S$  can be considered analogously.
  {\it Theorem $3.2$ is proved.}

\medskip
\medskip
{\it Remark.}
 We note that the following equality for Frech{\'e}t
derivatives is valid,
$$
d^n(RWR^{-1})(0)(z_1,...,z_n)
= Rd^n W(0)(R^{-1} z_1,...,R^{-1}z_n).
$$
We note also that due to holomorphity
(or commutativity with the operator of imaginary unit)
 Frech{\'e}t derivatives $d^n(RWR^{-1})(0)$
are $n$-\-linear forms with
respect to the field of complex numbers.



\medskip
\medskip
\medskip
\medskip
{\large 4.
Smoothed solution with small initial data
}
\medskip
\medskip

To prove that the maps $W,$ $S$ are holomorphic
for  large initial data we show that the positive (or
negative) frequency part of  solution of
the nonlinear equation smoothed
(in time) with
a test function with compact support in momentum space
is a polynomial.

\medskip
{\large Theorem 4.1.}
  {\it Let $\vartheta$ be from Theorem $3.2,$
$\vartheta>0$ (i.e. the Taylor expansion at zero
of the operators
 $RWR^{-1}$ and $RSR^{-1}$ converges for
  $\Vert z_{in} \Vert_{H^1(\R^3,\C)}<\vartheta).$
Let  $z_{in} \in H^1(\R^3,\C) \cap {\Sc}(\R^3,\C),$
$\Vert
z_{in}\Vert_{H^1(\R^3,\C)} < \vartheta,$
$ f \in {\cal F}(D(\R)),$ then
$$
\int dt f(t)RWU_0(t)R^{-1}z_{in}
$$
$$=\int dt
f(t)\sum^{N(f)}_{n=1}{1\over n!}
Rd^n W(0)(U_0(t)R^{-1}z_{in},...,
U_0(t)R^{-1}z_{in}),\eqno(4.1)
$$
where $N(f) = \max\{[p_0/\mbox{ {\it mass constant} }] + 1
\,\vert
\, p_0 \in \mbox{ supp }f^\sim\},$ $U_0(t)$
is the free dynamics (i.e. the dynamics defined by
the linear Klein-\-Gordon equation).
}

\medskip
\medskip
{\it Remark.}
  It is clear that $H^1(\R^3,\C)=D(\mu^{1/2})$ in
$H^{1/2}(\R^3,\C).$

\medskip
\medskip
{\it Proof of Theorem  $4.1.$}
  Since  $U_0(t)$ is an orthogonal (and symplectic)
transformation, so for
 $\Vert z_{in}\Vert_{H^1(\R^3,\C)} < \vartheta$
$$
\int dt f(t)RWU_0(t)R^{-1}z_{in}
$$ $$
= \int dt f(t)\sum^\infty_{n=1}{1\over n!}
Rd^n W(0)(U_0(t)R^{-1}z_{in},...,U_0(t)R^{-1}z_{in})
$$
$$
=\sum^\infty_{n=1}\int dt f(t){1\over n!}
Rd^n W(0)(U_0(t)R^{-1}z_{in},...,U_0(t)R^{-1}z_{in})
$$
(the series converges in $H^1(\R^3,\C)$ topology).

The
$n$-linear form
$$ n!^{-1} d^n RWR^{-1}(0)$$
is the transformation into  $H^1(\R^3,\C).$
The Schwartz theorem implies that
there exists a  unique generalized function
from the space ${\Sc}'(\R^{3n},\C)$ with values in
 $H^1(\R^3,\C)$
 (i.e. from the space ${\Sc}'(\R^{3n}, H^1(\R^3,\C)$)
 such that
$$ {1\over n!}Rd^n W(0)(R^{-1}z_1,...,R^{-1}z_n)
= R_n(z_1\otimes...\otimes
z_n) $$ for $z_1,...,z_n\in{\Sc}(\R^3,\C).$

Taking into account that
$U_0(t)R^{-1}z=R^{-1}(\exp(-i\mu t)
z),$ where
$\mu=(-\Delta+ m^2)^{1/2},$ we obtain $$ \int dt
f(t) {1\over n!} R d^n W(0)
(U_0(t)R^{-1}z_{in},...,U_0(t)R^{-1}z_{in}) $$
$$ =\int dt f(t)
R_n (\underbrace{\exp(-i\mu t)z_{in}\otimes...
\otimes\exp(-i\mu t)z_{in}}_{n}) $$ $$ =
R_n(f^\sim(\mu_1+...+\mu_n)z_{in}\otimes...
\otimes z_{in}), $$
here $\mu_j=(-\Delta_j+m^2)^{1/2}.$ The operator
$f^\sim(\mu_1+...+\mu_n)$
is a convolution operator,
 its Fourier transform (in spatial coordinates)
 is equal to $f^\sim(\sum_j
(p^2_j + m^2)^{1/2})$ and $ = 0$ for
$$n \geq N(f) = \max
\{[p_0/m] + 1 \, | \, p_0 \in
\mbox{ supp } f^\sim\}.$$
This implies that
 $$ \int dt f(t) RWU_0(t)R^{-1}z_{in} $$
$$= \int dt
f(t) \sum^{N(f)}_{n=1}{1\over n!}
Rd^n W(0) (R^{-1}\exp(-i\mu t)z_{in},...,
R^{-1}\exp(-i\mu t)z_{in}).\eqno(4.2) $$

The continuity allows to extend this equality
  on $$z_{in} \in H^1(\R^3,\C),\quad \Vert
z_{in}\Vert_{H^1(\R^3,\C)} < \vartheta.$$
Indeed, the real analyticity proved by Kumlin
\cite{Kumlin92}
and Theorem 26.2.4, ch. XXVI, \S1 (or Theorem
26.2.5, ch. XXVI, \S1, Theorem 26.2.6,
ch. XXVI, \S1) by Hille and Phillips,
see
\cite{HillePhillips}
 imply that all derivatives at zero are homogeneous
continuous polynomials, and, therefore, all polynomials
smoothed with  integration over
$f(t),$ are continuous also.
These results give the  continuity and Equality (4.1)
for
$\Vert
z_{in}\Vert_{H^1(\R^3,\C)} < \vartheta.$
  {\it Theorem $4.1$ is proved.}






\medskip
\medskip
\medskip
\medskip
{\large 5. Holomorphity on finite--dimensional
subspaces}
\medskip
\medskip

In this Section we consider the complex analyticity of
functions
$$
\int dt f(t)\langle
RWU_0(t)R^{-1}(z_{in}(\alpha)),h\rangle_{H^{1/2}(\R^3,\C)},\quad
f\in L_1(\R).
$$

\medskip
{\large Theorem 5.1.}
  {\it Let $z_{in}(\alpha)=\sum^N_{j=1}\alpha_j z_{in,j},$
$z_{in,j}\in H^1(\R^3,\C),$ $\alpha_j\in\C.$
Let $f \in
L_1(\R),$ $h \in H^{1/2}(\R^3,\C),$
then $$ \int dt f(t)\langle
RWU_0(t)R^{-1}(z_{in}(\alpha)),
h\rangle_{H^{1/2}(\R^3,\C)},
$$
is an entire antiholomorphic function in
$(\alpha_1,...,\alpha_N)\in\C^N.$
}

\medskip
{\it Proof of Theorem $5.1.$}
  For the proof we use Theorem 3.2 about complex
analyticity at the zero in
$H^1(\R^3,\C)$
($\simeq D(\mu^{1/2})$ in  $H^{1/2}(\R^3,\C)$),
Equality
(4.1)
 and the real analyticity proved by
Kumlin
\cite{Kumlin92}.

Let $f\in{\cal F}(D(\R))$ and
\be
F_1(z_{in}(\alpha),f)
&=&
\int dt f(t)\langle
RWU_0(t)R^{-1}(z_{in}(\alpha)),h\rangle_{{\cal H}_S},
\\
F_2(z_{in}(\alpha),f)
&=&
\int dt f(t) \sum^{N(f)}_{n=1}
{1\over n!}\langle
Rd^nW(0)(U_0(t)R^{-1}z_{in}(\alpha),...,
U_0(t)R^{-1}z_{in}(\alpha),h
\rangle_{{\cal H}_S}.
\ee
Then  $F_1(z_{in}(\alpha),f)$ is
a real analytic function for all
 $\alpha\in\C^N$ and $F_2(z_{in}(\alpha),f)$ is an
entire antiholomorphic function in
 $\alpha\in\C^N.$ The real analyticity of
 $F_1(z_{in}(\alpha),f)$ follows from the proof of real
analyticity given by
Kumlin
\cite{Kumlin92},
from  energy
conservation and  uniform in $\alpha$
approximation  of the integral over time
by finite sums. The function
$F_2(z_{in}(\alpha),f)$ is an entire antiholomorphic
function. This is implied by  the facts that this
function is a polynomial, generated by
Frech{\'e}t derivatives at zero  and by holomorphity
of the transformation  $RWR^{-1}$ in a zero neighborhood
in  $H^1(\R^3,\C),$ see  Theorem 3.2.

Theorem  4.1 implies that there exists $\vartheta(f)>0$
(and
depending on $z_{in,1},...,z_{in,N})$ such, that
$$
F_1(z_{in}(\alpha),f)=F_2(z_{in}(\alpha),f)
$$
for $|\alpha|<\vartheta(f).$ Then Lemma 5.2 about
 uniqueness implies the equality
$$
F_1(z_{in}(\alpha),f)=F_2(z_{in}(\alpha),f)
$$
for all $\alpha\in\C^N$ and, therefore
the function $$
\int dt f(t)\langle
RWU_0(t)R^{-1}(z_{in}(\alpha)),h\rangle_{H^{1/2}(\R^3,\C)}
$$
is an entire antiholomorphic function in $\alpha.$

If now $f\in L_1(\R),$ then there exists
a sequence  $f_n\in{\cal F}(D(\R)),$ $f_n\to f$ in
$L_1(\R)$ and
$$ F_1(z_{in}(\alpha),f)=\lim_n
F_1(z_{in}(\alpha),f_n).\eqno(5.1) $$
 Due to the estimate $$
|F_1(z_{in}(\alpha),f-f_n)|\leq\int dt|f(t)-f_n| \,
\Vert z_{in}(\alpha)\Vert_{H^1(\R^3,\C)} \,
\Vert\mu^{-1/2}h\Vert_{H^1(\R^3,\C)}
$$
 uniform in $\alpha$ for bounded  $\alpha,$
the convergence in (5.1) is uniform in $\alpha$ and,
 thus, $F_1(z_{in}(\alpha),f),$
 $f\in L_1(\R),$ is an entire
 antiholomorphic function in $\alpha.$
  {\it Theorem $5.1$ is proved.}

\medskip
\medskip
Now we prove the statement about uniqueness. We
prove this statement in the following form.

\medskip
\medskip
{\large Lemma 5.2.}
  {\it Let $A(\alpha)$ and $B(\alpha)$ be two real analytic
functions from  $\R^N$ into a Banach space ${\bf B}.$ Let for
 $|\alpha| < a$ with some $a>0$ $A(\alpha) = B(\alpha).$ Then
$A(\alpha) = B(\alpha) $ for all $\alpha \in \R^N.$
}

\medskip
\medskip
{\it Remarks.}
  1. Real analyticity means a convergence at any point of
local  Taylor
expansion with Frech{\'e}t derivatives.

2. For us it is sufficient
to consider the case  ${\bf B}=\C.$

\medskip
\medskip
{\it Proof of Lemma $5.2.$}
  Let $ Anal(\alpha,r)$ be an (open) ball with center at
$\alpha$ and radius $r.$ Let $r_{A,B}(\alpha): =
\min(r_A(\alpha), r_B(\alpha)),$
where $r_A(\alpha), r_B(\alpha)$ are
radii of convergence
of the Taylor expansion with the center at
$\alpha$ for the function $A$ and
for the function  $B,$
 respectively.
By  condition of Lemma 5.2
 $r_{A,B}(\alpha)>0$ for all  $\alpha\in\R^N.$
 Let ${\cal O}(0,r)$
be an (open) ball with radius $r$ having the properties
 $$ A(\alpha)=B(\alpha)\quad\mbox{ for }\quad
|\alpha|<r.  $$ By conditions of Lemma 5.2
$A(\alpha)=B(\alpha)$ in $$Anal(0,r_{A,B}(0)) \cap
\{\alpha\,\vert\,|\alpha|<\vartheta, \alpha \in \R^N\}.$$
Thus,  there exists a ball
${\cal O}(0,r)$ with radius $r>0.$

Let $$r_{max}=\sup\{r\;\vert\;
	\exists \;{\cal O}(0,r), \, A(\alpha)=B(\alpha)
	\;\forall\alpha\in {\cal O}(0,r)\}.$$
We show that $r_{max}=\infty.$

It is clear that
$r_{max}>0.$ If $r_{max} <\infty,$ then continuity
and the real analyticity imply
that
 $ A(\alpha) = B(\alpha)$ and
$\partial^j A(\alpha) = \partial^j B(\alpha)$ for
 $\alpha \in \overline{{\cal O}(0,r_{max})}$
  and all $j.$
 The compact $$ \{\alpha \in \R^N \;
 \vert\; |\alpha| = r_{max}\} \subset \bigcup_{|\alpha|=r_{max}}
 Anal(\alpha,{1\over 4}r_{A,B}(\alpha)) $$
and compactness of
  $\{\alpha\in\R^N \;\vert \;
	|\alpha|=r_{max}\}$
 implies
that there exists a finite number
of balls
 $Anal(\alpha_k, {1\over 4} r_{A,B}
 (\alpha_k)),$ $k=1,...,K,$ such that  $$ \{\alpha \in
\R^N \; \vert \; |\alpha| = r_{max}\} \subset
\bigcup_{1\leq k\leq K}
Anal(\alpha_k, {1\over 4} r_{A,B}(\alpha_k)).  $$

Let
$${\cal O}(K) =
{\cal O}(0,r_{max})\cup \bigcup_{1\leq k\leq K}
 Anal(\alpha_k, {1\over 4} r_{A,B}(\alpha_k)).$$
 Let
 $$r(K) =
 {1\over 4}\min_{1\leq k\leq K} r_{A,B}(\alpha_k).$$
 Then
 $r(K)>0$ and ${\cal O}(0,r_{max}+ r(K)) \subset{\cal O}(K).$
Really, if $\alpha \in
 {\cal O}(0,r_{max}+r(K)),$
 $|\alpha| \geq
r_{max},$ then the vector
 $${\alpha \over |\alpha|} r_{max} \in
 Anal(\alpha_k, {1\over 4}r_{A,B} (\alpha_k))\cap {\cal O}
 (0,r_{max})$$ for some  $k.$
 Then $\alpha = \alpha_k +
\beta$ and
\be |\beta| &\leq&
|\alpha-{\alpha\over|\alpha|} r_{max}|+|{\alpha\over
|\alpha|}r_{max}-\alpha_k|
\cr
&\leq& |\alpha|-r_{max}+|{\alpha
\over|\alpha|}r_{max}-\alpha_k|
\cr
&\leq& r(K)+{1\over 4}r_{A,B}(\alpha_k)
\cr
&\leq&{1\over 2}r_{A,B}(\alpha_k).
\ee
Therefore,
$$
A(\alpha)=\sum{1\over j!} \partial^j  A(\alpha_k)
(\beta-\alpha_k)^j
=\sum{1\over j!} \partial^j
B(\alpha_k)(\beta-\alpha_k)^j=B(\alpha), $$
because the series converges in
$ Anal(\alpha_k, r_{A,B}(\alpha))$
and all derivatives for
functions $A$ and $B$ at point $\alpha_k$ coincide.
{\it Lemma $5.2$ is proved.}

\medskip
\medskip
{\large Lemma 5.3 (Consequence of Lemma 5.2).}
  {\it Let  $A(\alpha)$ be a complex holomorphic
in $\C^N$ function with
values in a complex Banach space ${\bf B}$ and let
$B(\alpha)$ be a real analytic function in $\C^N$ with
values in the same complex Banach space and
$A(\alpha) = B(\alpha)$ for $|\alpha|<a$
for some  $a>0.$
Then  $A(\alpha) = B(\alpha)$ for all $\alpha\in\C^N$
and $B(\alpha)$ is a complex holomorphic function.
}

\medskip
\medskip
{\it Remark.}
  By the space $\C^N$ we mean the space with standard basis,
given, for instance, in the form
$(z_1,...,z_N)\in \C^N,$
$z_j= x_j+iy_j,$ $x_j, y_j \in\R.$
 The complex analyticity is
holomorphity with respect to the variables
 $(z_1,...,z_N)$ and the real analyticity
is an expansion into a (local) Taylor series
with respect variables
 $(x_1, y_1,..., x_N, y_N).$



\medskip
\medskip
\medskip
\medskip
{\large 6. Unitarizability of operators
$dW(u_{in})$ and $dS(u_{in})$}
\medskip
\medskip

In this Section we complete the proof of the uniqueness of
essential unitarizability of operators
 $dW(u_{in})$ and $dS(u_{in})$  for small solutions, i.e.
Theorem 3.1.
For this purpose we
give here more detailed proof of
required assertions. These assertions are Theorem
6.1, 6.2, 6.3
(Theorem 6.3 is closely connected with Theorem 3.1)
and
these theorems have been used by
Paneitz
\cite{Pan80-82,Pan81,Pan82}
for the proof of an analogue of Theorem 3.1.

\medskip
\medskip
{\large Theorem 6.1} (see
\cite[ch. III, \S 1, \S 2]{DaletskiiKrein71}
and also
Theorem 6
\cite{Pan81}).
  {\it Let $t\to A(t)$ be a strongly continuous
norm-\-bounded map from $\R$ to ${\cal L}({\cal H})$
($=$ the algebra of bounded operators on
(the real Hilbert space) ${\cal H}$) such that
$$
N(A):=\int^{+\infty}_{-\infty}\Vert A(t)\Vert dt<2.
$$
Then
$$
U(t;\lambda)=I+\lambda\int^t_{-\infty} A(s) U(s;\lambda)ds
$$
has a unique continuous solution for all
$\lambda;$
$\Vert U(t;\lambda) \Vert
\leq \exp (|\lambda| N(A))$
for  all $t;$
$W(\lambda) := U(0;\lambda),$
$\lim_{t\to+\infty} U(t;\lambda)$
exists in norms,
$S(\lambda):= \lim_{t\to+\infty} U(t;\lambda).$
If $|\lambda|<2/N(A)$ then
$(W(\lambda) + I)^{-1}$ and
$(S(\lambda) + I)^{-1}$ exist,
$\lambda \to (W(\lambda) + I)^{-1}$
and
$\lambda \to (S(\lambda) +I)^{-1}$ are
analytic, and
$$ (W(\lambda)- I)(W(\lambda) + I)^{-1}
$$ $$= 2\int^\lambda_0
\int^0_{-\infty} (I+W(\rho)^{-1})^{-1}
U(s;\rho)^{-1} A(s)
U(s;\rho)(I+W(\rho))^{-1}dsd\rho,\eqno(6.1) $$
$$
(S(\lambda)- I)(S(\lambda) + I)^{-1}
$$ $$= 2\int^\lambda_0
\int^{+\infty}_{-\infty} (I + S(\rho)^{-1})^{-1}
U(s;\rho)^{-1} A(s)
U(s;\rho) (I + S(\rho))^{-1} dsd\rho.\eqno(6.2) $$
}

\medskip
\medskip
{\it Remark.}
  In this chapter notations $W(\lambda),$ $S(\lambda)$
denote linear operators and correspond to the first
derivatives of (nonlinear) operators $W(u_{in}),$
$S(u_{in})$ which are wave and scattering operators
of the nonlinear equation (1.1).

\medskip
\medskip
{\it Proof of Theorem $6.1.$}
  The first part is well known,
and follows immediatately from the
usual ``ti\-me-\-or\-der\-ed
exp\-on\-en\-ti\-al" \, form of the
solution $U(t;\lambda),$ a norm-\-con\-verg\-ent
power series in
$\lambda,$ we refer to the
Krein proof, see
\cite{DaletskiiKrein71}.

In the second part we follow to
Paneitz
\cite[Theorem 6]{Pan81}.
and give  the proof for the case
of the operator $W(\lambda).$

We  need some notation. Define
$$
h(t)=\cases{+{1\over 2} &for  $ t\geq 0,$\cr
-{1\over2} &for $t<0.$  \cr}
$$
This is the Green function of
the linear differential
equation of first order in time.

Let ${\cal B}$  be the Banach space
of continuous functions $f(t)$ from
$\R$ to ${\cal H}$
having continuous limits
$f(\pm\infty)$ as $t\to\pm
\infty,$ with the sup norm
$$
\Vert\vert f\Vert\vert=\sup_{t\in\R}\;
\Vert f(t)\Vert.
$$
We define an integral operator
$K_-:{\cal B}\to{\cal B}$ by
$$
(K_-f)(t)=\lambda\int^0_{-\infty} h(t-s) A(s) f(s)ds.
$$
If $ N(A) < 2$
and $|\lambda|<2/N(A)$ clearly $\Vert\vert
K_-\Vert\vert<1,$ and
then
\footnote{For the case of the wave operator
it is sufficient to take the condition
$N_- := \int^0_{-\infty}\Vert A(t)\Vert dt <
2$.}
$$ (I-K_-)g=f\eqno(6.3) $$
(for a given  {\it constant}
$f\in{\cal H}$)
has a unique solution $g\in{\cal B}.$
This solution $g$ satisfies the boundary conditions
$$
g(-\infty)+g(0)=2f\eqno(6.4)
$$
and the equation
$$
g(t)=g(-\infty)+\lambda\int^t_{-\infty}
A(s)g(s)ds.\eqno(6.5)
$$
The boundary condition (6.4)
is implied by Equation (6.3)
and gives equalities
$$ g(0)-(K_-g)(0)=g(0) -
{1\over 2}\lambda\int^0_{-\infty} A(s)g(s)ds = f,
$$ $$
g(-\infty) - \lim_{T\to -\infty}(K_-g)(T)
=g(-\infty)+{1\over2}\lambda
\int^0_{-\infty}A(s)g(s)ds = f.
$$
In the continuous case Equation (6.5) is
equivalent to the differential equation.
 Equation (6.5) can be obtained from
Equation (6.3)
by differentiation in $t.$

Then by uniqueness $g(t) =
U(t;\lambda) g(-\infty),$
$g(0) = W(\lambda)g(-\infty),$
and $2f = (I + W(\lambda)) g(-\infty).$
$g$, hence $g(-\infty),$
depends continuously on $f,$ so
$\lambda\to (I+W(\lambda))^{-1}$
is an analytic map
into ${\cal L}({\cal H})$ for
$0 \leq \lambda \leq 1,$
that is, a power
series convergent in norm,
this is implied by the inequality
$\Vert\vert K_-\Vert\vert<1$ for
$|\lambda| < 2/N(A).$
We emphasize that ${\cal H}$
is a real Hilbert space
(see
\cite[p. 316]{Pan81},
thus $\lambda$ is a real constant also.

It is easily checked that
$$
{d\over d\lambda} U(t;\lambda)
=U(t;\lambda) \int^t_{-\infty}
U(s;\lambda)^{-1} A(s) U(s;\lambda)ds
$$
(by solving the first order differential equation
which the l.h.s. satisfies),
 $$ {dU(t;\lambda)\over d\lambda}=\int^t_{-\infty}
A(s) U(s;\lambda)ds
+ \int^t_{-\infty} A(s){dU(s;\lambda)\over
d\lambda}ds, $$

$$
{d^2 U(t;\lambda)\over dt d\lambda}
= A(t) U(t;\lambda) + \lambda A(t)
{dU(s;\lambda)\over d\lambda},
$$
see
\cite[ch. III, \S 1, s. 4, (1.1),
(1.10), (1.19)]{DaletskiiKrein71}.
This equation corresponds to the linear equation with an
external field.

Thus, $dW(\lambda)/d\lambda =
W(\lambda)R_-(\lambda),$
where
$$
R_-(\lambda) = \int^0_{-\infty}
U(s;\lambda)^{-1} A(s) U(s;\lambda)ds.
$$
Defining $X_-(\lambda) =
(W(\lambda)-I)
(W(\lambda)+I)^{-1},$ clearly
$$
{d\over d\lambda}X_-(\lambda)
= 2 (I+W(\lambda))^{-1}
\, W(\lambda) \, R_-(\lambda)(I+W(\lambda))^{-1},
$$
and since $W(0)=I,$ $X_-(\lambda)$
has the integral expression
(6.1).
  {\it Theorem $6.1$ is proved.}

\medskip
\medskip
In the case when $A(t)\in sp({\cal H})$ it is well known that
then $U(t;\lambda),$ $W(\lambda),$ $S(\lambda )\in Sp({\cal H}).$
 Here we use Paneitz's notation \cite{Pan81},
 $Sp({\cal H})$
 $(sp({\cal H}))$ is the group (respectively, Lie
 algebra) of bounded invertible
  (respectively, bounded) operators on
${\cal H}$ preserving (respectively, skew with
respect to) the symplectic form $a(.,.).$

\medskip
\medskip
{\large Theorem 6.2.}
(see
Theorem 7
\cite{Pan81}).
  {\it Let $t\to A(t)$  be a strongly
continuous norm-bounded map from
$\R$ into $sp({\cal H})$
taking values in the closed positive
cone $\overline{C_0}$
(the notation follows to
\cite[pp. 316, 319]{Pan81}).

($1_W$)
If $\int^0_{-\infty}\Vert A(t)\Vert dt<2$ then the wave
operator $W:$ $X(-\infty) \to X(0)$ for the equation
$dX/dt = A(t)X$
is the Cayley transform of a $Y_- \in \overline{C_0}.$

($1_S$) If $\int^{+\infty}_{-\infty}
\Vert A(t)\Vert dt<2$ then
the scattering operator $S:$
$X(-\infty) \to X(+\infty)$  for the
equation $dX/dt = A(t)X$ is the Cayley transform of a
$Y \in \overline{C_0}.$

($2_W$) If in addition $\cap_{{t\in\R, t\leq 0}}\{\ker
A(t)\} = \{0\},$ then also $a(Y_- v,v) > 0$
for all $v\not=0$
(in Paneitz's notations
\cite{Pan80-82}).

($2_S$) If in addition
$\cap_{t\in\R}\{\ker A(t)\} = \{0\},$
then also $a(Yv,v) > 0$ for all $v \not= 0$
(in Paneitz's notations
\cite{Pan80-82}).

($3_W$) If furthermore $\int^0_{-\infty} A(t) dt$
(as a strong-operator topology integral)
is in the interior $C_0,$ then
$Y_-\in C_0$ also, and then
$W$ commutes with a unique complex
structure, i.e., is uniquely unitarizable.

($3_S$) If furthermore
$\int^{+\infty}_{-\infty} A(t) dt$ (as a
strong-operator topology integral) is in
the interior $C_0,$ then $Y \in C_0$ also, and
then $S$ commutes with a unique complex structure, i.e.,
is uniquely unitarizable.
}

\medskip
\medskip
{\it Remarks.}
  1. The conditions
 $\int^0_{-\infty}\Vert A(t)\Vert dt<2$ and
 $\int_0^{+\infty}\Vert A(t)\Vert dt<2$ are sufficient for the
consideration of wave operators. In the quantum case
this is connected with an elastic scattering for energies less
than $4m.$ We intend to consider this question later
in connection with the unitarity of the quantum scattering
operator.

2. The restriction on the subspace of kernels in the
condition of Theorem 6.2 corresponds to the exlusion of ``bound
states", i.e. vectors $v\not= 0$ such that $Sv=v.$

\medskip
\medskip
{\it Proof of Theorem $6.2.$}
  ($1_W$) It follows from the integral formula in
 Theorem 6.1 for
 $(W-I)(W+I)^{-1},$
 the invariance of $\overline{C_0}$ under
$(U_W)^{-1}:$
 $T\to 2(I+W^{-1})^{-1} T(I+W)^{-1}$
(Theorem 5
\cite{Pan81})
and $Ad \,W$ for
 $W \in Sp({\cal H}),$ and the fact
 that $\overline{C_0}$ is a
 convex cone.

($2_W$) For $\rho =0$ the integrand in
(6.1)  is equal to
${1\over 2}\int^0_{-\infty} A(t)dt.$
Thus, if for some
$v\in{\cal H}\,$ $\int^0_{-\infty}
a(A(t)v,v)=0,$ then
$a(A(t)v,v) = 0$ for all $t\leq 0,$
due to continuity in $t$ and the inclusion
$A(t) \in \overline{C_0}.$
 Since
$A(t) \in sp({\cal H}),$ so $a(A(t)v,v) = {\Sc}_0(J^{-1}
A(t)v,v)$
(in the notations of
\cite[pp. 316-319]{Pan81})
and $J^{-1} A(t)$ is a symmetrical and positive
operator.
Therefore,  the equality
$a(A(t)v,v) = 0 =
{\Sc}_0(J^{-1} A(t)v,v)$
 means that $(J^{-1}
A(t))^{1/2} v = 0,$ and,
thus,
 $J^{-1} A(t)v = 0$ and
$A(t)v = 0,$
that is,
 $v \in \cap_{t\leq 0} \{\ker A(t)\} =
\{0\}$
and $v=0.$
This means that if $v \not= 0,$ then, due to
(6.1), $a(Y_- v,v)>0.$

($3_W$) As noted before, the inner integrand in (6.1) is
norm-continuous in $\rho.$
Thus the integral is in the interior
$C_0$ if
$\int^{0}_{-\infty} A(t) dt \in C_0,$ by the
observation in part ($2_W$).
Finally, apply
Theorem 2
\cite[Theorem 2, p. 318]{Pan81}
 to
$(Cayley\;\; transform)^{-1}(W).$

 The consideration of ($1_S$),  ($2_S$), ($3_S$),
i.e. the
integral formula  (6.2) in Theorem 6.1 for
 $(S-I)(S+I)^{-1}$
is completely the same as the case of (6.1).
  {\it Theorem  $6.2$ is proved.}

\medskip
\medskip
{\large Theorem  6.3}
($=$ Theorem 4.4B
\cite{Se66},
\cite[p. 491]{Se68},
\cite[p. 110]{Pan80-82}).
  {\it Let $u_{in}$ be a given finite--energy
solution of the free equation
$$ \Box u_{in}+m^2 u_{in}=0 \eqno(6.6)
$$
such that $grad \,u_{in}$
is also of finite energy, and suppose
that the
$L_\infty$-\-norms (over space) of $u_{in}(t,\cdot)$
and $grad \,u_{in}(t,\cdot)$
are bounded by $ const \,(1 + |t|)^{-3/2}.$
Then if $u_{in}$
 in a certain norm
 is sufficiently small,
 there
exist unique solutions $u$ and
$u_{in},$ $u_{out}$ of
$$ \Box u+m^2 u + \lambda u^3 =0
$$ and $(6.6)$
respectively, such that
$$ \Vert u-u_{in}\Vert_{H^1\oplus L_2}
\to 0\quad\mbox{
as } \quad t\to -\infty, $$
$$ \Vert u-u_{out}\Vert_{H^1\oplus L_2}
\to 0\quad\mbox{
as } \quad t\to +\infty, $$
and
$ u,$ $u_{in},$ $u_{out} = O(|t|^{-3/2})$
in $ L_\infty(\R^3). $
}

\medskip
 It may also be  discerned from these method that in
these and other similar situations
$$ |t|^{3/2} \,
\Vert u - u_{in}\Vert_\infty \to 0
\quad\mbox{ as } \quad t\to -\infty,
$$
 $$ |t|^{3/2} \,
\Vert u - u_{out}\Vert_\infty \to 0
\quad\mbox{ as } \quad t\to +\infty,
\eqno(6.7)
$$
a fact we will need later.

\medskip
\medskip
{\large Corollary 6.4} (See also Theorem 4.4B
\cite[p. 110]{Pan80-82},
\cite[p. 491]{Se68},
Corollary
\cite[p. 115]{Pan80-82}).
  {\it Take nonvanishing solutions
$u,$ $ u_{in},$ $ u_{out}$ as in
Theorem $6.3$
earlier (which implies in particular
that $u$ is uniformly bounded and
$\int^{+\infty}_{-\infty} \Vert
u (t, \cdot)^2\Vert_\infty dt
< \infty$), and assume that
$\int^{+\infty}_{-\infty}
\Vert u (t, \cdot)^2\Vert_\infty dt < 2m$
and that
for some time $t$
$$
(m^2-\Delta)^{5/4} u_{in}(t,\cdot),
\quad (m^2-\Delta)^{3/4} u_{in}(t,\cdot)
\in L_1(\R^3),
\eqno(6.8)
$$
then $dW(u_{in}) :
(H^{1/2}\oplus H^{-1/2},J)
\to (H^{1/2} \oplus H^{-1/2},J)
$
and
$dS(u_{in}) :
(H^{1/2}\oplus H^{-1/2},J)
\to (H^{1/2} \oplus H^{-1/2},J)
$
are uniquely essentially unitarizable.
}

\medskip
\medskip
{\it Proof
of Corollary $6.4.$}
  To prove Corollary we use   Theorem 2b)
\cite[p. 115]{Pan80-82})
(= Theorem 7(2)
\cite{Pan81}
= Theorem 6.2
\cite{Pan82},
 see also
Theorem 6.3
\cite{Pan82}),
 it remains only to show
that all
$$
-A(t) = e^{-tQ} \pmatrix{0&0\cr F^\prime (u)&0\cr} e^{tQ},
\quad
\quad
Q = \pmatrix{0&1\cr -\Delta+m^2&0\cr},
$$
vanish on no nonzero vector $(v_1,v_2) \in
H^{1/2}\oplus H^{-1/2}.$
Now $e^{tQ} (v_1, v_2) = (v(t), {\dot v}(t))$
for
$(v_1,  v_2)$ with finite
$H^{1/2}\oplus H^{-1/2}$
norm satisfying the Klein-\-Gordon equation, and
$F^\prime(u) v : = u^2 v = 0$
is equivalent to $u v = 0.$
The conclusion $v \equiv 0$ follows from
(6.7),
the asymptotics of $u_{in}(t,0)$
in some Lorentz frame
\cite[Corollary 2]{SNelson71}
obtainable from
(6.8),
hyperbolicity, and
the vanishing of any solution $v,$
$(v(t), {\dot v}(t)) \in H^{1/2} \oplus H^{-1/2},$
which vanishes in a backward ti\-me-\-li\-ke cone,
see Theorem 7.1
 (see also
\cite{PanSe80,Goodman64},
\cite[Corollary $1^\prime$]{SNelson71}.
 Our proof  uses also
 Theorem 6.3
\cite[Theorem 6.3]{Pan82},
see also  Theorem 2
\cite[Theorem 2, p. 115]{Pan80-82}.

To consider the case (of the derivative)
of  wave operator
$dW(u)$ we apply
 Theorem 7,
\cite[Theorem 7]{Pan81},
for the interval $(-\infty,0]$ (or $[0,\infty)$ for
the forward  wave operator)\footnote{It is possible
to  use  the interval
 $(-\infty,T],$ or $[T,\infty),$ also.}
and use Theorem 6.3
\cite[Theorem 6.3]{Pan82}
and Theorem 2
\cite[Theorem 2, p. 115]{Pan80-82}.

In Theorem 7.1 we give the detailed  proof
of the statement
that the condition
$$u(t,x) v(t,x) = 0$$ for all
$(t,x)\in\R_-\times\R^3$ implies that  $v(t,x) = 0$ for all
$(t,x) $ belonging to
some backward time--like cone.
The last condition implies that $v = 0.$
The ideas of the proof of Theorem  7.1
 are analogous to
\cite{Pan81,Pan82,PanSe80,Pan80-82}
and use the S.~Nelson results
about asymptotic behavior of free solutions,
see
\cite[Corollary 2, Corollary $1^\prime$]{SNelson71}.
   {\it Corollary $6.4$ is proved.}

\medskip
{\it Remarks.}
  1. We point out  that the paper
\cite{PanSe80}
of Paneitz and Segal formulates a statement, gives
 some references, but does not contain the
proof of this statement.
The  more detailed exposition of Paneitz
\cite{Pan80-82}
formulates the assertion $ v = 0$
in the forward cone
(for the case of  derivative of the
scattering operator) as a condition.
This is not quite appropriate because
only  the condition
$v(t,x) = 0$ on the forward time-\-like cone appears
naturally
 \cite{PanSe80}.
The last condition implies
that  $v_{in}(t,x) = 0$ on the forward cone and,
therefore,
$v_{in} = 0$
\cite{Mor63}
and $v = dS(u_{in}) v_{in} = 0$
(and also by Goodman
\cite{Goodman64}
$v(t,x)$ ($=(dS(u_{in}) v_{in})(t,x)$
is equal to zero on the cone).
However the Paneitz's paper
\cite{Pan80-82}
uses the absence of solution
with initial data from
$H^{1/2} \oplus H^{-1/2}$
and
 the Morawetz's paper
\cite{Mor63}
contains the proof for solution with finite energy only.
In this case the use of Theorem 1
\cite[Theorem 1]{SNelson71}
is more appropriate.
Theorem 1
\cite[Theorem 1]{SNelson71}
 gives the required assertion about the statement
that a solution equal to zero in some
 (backward) light
cone is the zero solution.

2. The part  of Corollary 6.4 is contained in
Theorem 16.3
\cite{Pan82}.

3. We remark that to use
Corollary
\cite{Pan80-82}
 we
take the space $R(H^1\oplus L)$
as the subspace $D$ ($= D(\mu^{1/2})$ in
$H^{1/2}(\R^3,\C)$)
(see Definitions
\cite{Pan80-82}
and also
\cite{Pan82}).

4. To apply Theorem 6.3
\cite{Pan80-82}
we use that the condition
$a(Sv,v) >  0$  for all $v \not= 0$
is implied by Condition 2 of Theorem 7
\cite{Pan81}
(i.e. $\cap_{t\in\R}\{\ker A(t)\} = 0$)
and by Statement 2 of Theorem 7
\cite{Pan81}
(for the wave operator this condition
is
$\cap_{t\in{\R}, t\leq 0}\{\ker A(t)\}=\{0\}$
and $a(Wv,v) > 0$ for all $v \not = 0$).
 Condition 7(2) ($=$
Statement  7(2)
\cite{Pan81})
and the simple
algebraic identity
$$ a((S - I)(S + I)^{-1}v,v) = 2a(Sw,w), $$
where $w = (S + I)^{-1}v,$
imply positivity of the inverse Cayley  transform
$$ Y = (S - I)(S + I)^{-1},\quad\;\;
Y_- = (W_{-} - I) (W_- + I)^{-1}.$$
Here (and in
\cite{Pan82})
the  symplectic form
is denoted by $a(\cdot,\cdot)$ instead of
the notation $\omega(\cdot,\cdot)$
that is used in the other places.

5. For the consideration of d-dimensional case
 we intend to give the detailed
(and independent) proof of
 Theorem 6.3
\cite{Pan82}
 and Theorem 2
\cite{Pan80-82}
  separately.




\medskip
\medskip
\medskip
\medskip
{\large 7. $\bigcap\{\ker A(t)\}=\{0\}.$
A solution equal zero
in an infinite time column   is equal zero}
\medskip
\medskip

To complete the proof of assertions about holomorphity
we need the following statement.

Let us consider the system of equations
$$
\left.\begin{array}{cc}
\Box u+m^2 u+\lambda u^3&=0\\
\Box v+m^2 v+3\lambda u^2 v&=0
\end{array}\quad \right\}\eqno(7.1)
$$
(solutions of with are the tangent space
of the manifold of solutions of the nonlinear
equation).

\medskip
{\large Theorem 7.1.}
 {\it Let $(u,v)$ be solutions of the system of
equations $(7.1).$
and let for the initial $in$-data of the solution
$u$ satisfies hypotheses of Theorem $3.1.$
Let $u \not = 0$ and $u(t,x)v(t,x) = 0$
for all $(t,x)\in\R_-\times\R^3.$
Then $v=0.$

The analogous statement is valid for
the scattering operator
and for the $out$--wave operator.
}

\medskip
{\it Remark.}
  Here, as in Theorem  3.1,
 it is sufficient to take initial data with weaker
conditions than in Theorem 3.1, see
Paneitz
\cite[Corollary, p.115-116]{Pan80-82},
 Paneitz, Segal
\cite{PanSe80}.

These conditions are
 required for the consideration of
solutions behavior. If the space
$\Sc_{Re}(\R^3) \oplus \Sc_{Re}(\R^3)$
is used as the space of
 initial $in$-\-data these
conditions  are  fulfilled. In any case
these conditions are fulfilled for sufficiently small initial
 $in$-\-data.

\medskip
{\it Proof of Theorem $7.1$} (the idea of the proof is given in
\cite{Pan80-82}).
 In our case, i.e.  for Equation
(7.1),
the operator $A(t)$ in the statement
$$
\bigcap_{\scriptstyle t\in\R,\atop
\scriptstyle t\leq 0}\{\ker A(t)\}=\{0\}
$$
(see
\cite{Pan80-82,Pan81,Pan82}
 and the Paneitz's notation)
corresponds to the operator
$$
U_0(-t) \pmatrix{0&0\cr -3\lambda u(t)^2&0\cr} U_0(t).
$$

This correspondence fulfills in the interaction
representation for  the second equation in (7.1).

Therefore, the condition  $A(t)v(t) = 0$ is reduced
(for the $in$-wave operator)
to the condition
$$
u(t,x)v(t,x)=0\;
\mbox{ for all }\; (t,x)\in\R_-\times \R^3.
$$
Analogously, for the $out$-wave operator the condition
$$
u(t,x)v(t,x)=0\;\mbox{ for all }\; (t,x)\in\R_+\times \R^3
$$
appears and for the scattering operator the condition
$$
u(t,x)v(t,x)=0\;\mbox{ for all }\; (t,x)\in\R\times \R^3
$$
appears.

The proof of the fulfillment of the condition
$$
\bigcap_{\scriptstyle t\in\R \atop \scriptstyle t\leq 0}
\{v \,|\, u(t,x)v(t,x) = 0 \quad\; \forall x\}=\{0\},
$$
or
$$
\bigcap_{\scriptstyle t\in\R \atop \scriptstyle t\geq 0}
\{v \,|\, u(t,x)v(t,x) = 0 \quad\; \forall x\}=\{0\},
$$
follows from hyperbolicity, relations
$$
|t|^{3/2}\,\Vert u-u_{in}\Vert_\infty\to 0\;
\mbox{\quad for \quad}\; t\to-\infty,
\eqno(7.2)
$$
or
$$
|t|^{3/2}\,\Vert u-u_{out}\Vert_\infty\to 0\;
\mbox{\quad for \quad}\; t\to+\infty,
\eqno(7.3)
$$
and the behavior of $u_{in}(t,x),$ or $u_{out}(t,x),$
for large times, see
\cite[Corollary 2]{SNelson71}
 and Theorem 7.2.

For this purpose, we prove the existence
of some time column in which $v(t,x)$ are equal zero.
To proof this fact we use hyperbolicity and mentioned
convergences.
This time column, i. e. time-\-like positive axes with
vertices on some space-\-like base,
is directed to the
past for the $in$-wave operator (and to the future
for the $out$-wave operator).
The base of this column is arranged sufficiently
far in the past.
Since  $v(t,x)$ for $t\leq 0$  satisfies the
free equation and the causal envelope
of this column is the backward light cone, so
the solution is equal to zero in
 the backward light cone and, thus, it is
the zero solution
\cite[Corollary 1]{SNelson71}.

This column is defined by the solution
 $u$ and can be constructed with the help
of the  solution $u_{in}$, namely, with
 the help of (neighborhood of) the point in
momentum space at which
$(u^+_{in})^\sim(k) \not = 0$
(we note that due to the condition
(6.8)
 $(u^+_{in})^\sim(k)$ is
continuous in $k$).
With the help
of the S.~Nelson's result
\cite[Corollary 2]{SNelson71},
which describes the behavior of a free solution
for large times and approximating the solution
 $u = W(u_{in})$
by its  $in$-data we obtain that
 $u(t,x) \not = 0$ for
$(t,x)$ belonging to some set. This set, after
the certain  Lorentz rotation has the form
  $\Sigma (a,r) =
 \cup [t_n,t_n +a]
 \times
 \{x\in\R^3 \,\vert\, |x| \leq r\}.$
 Therefore, $v(t,x) = 0$ on $\Sigma (a,r).$
 Since for
 $t \le 0$ $v(t,x)$
 satisfies the free equation, hyperbolicity
 implies that the parameters
$a$ and $r$ can be taken such that the causal envelope
$\Sigma (a,r)$ contains some
(more narrow than $r$) column,
and, thus, it is the backward light cone
\cite[ch. 5, \S 28]{Vla66}.

\medskip
{\it Remarks.}
 1. We remark that the time column (to the past) is
an infinite time cylinder of the form
 $(a_0, {\bf a}) + \R_+ \times \{x \in \R^3 \, |\,
|x| \le r\}.$

2. The causal envelope for the nonlinear wave equation,
 or for the system of equations  (7.1),
is  the light come  envelope as
in the case of the linear wave equation. This form of
   causal envelope
 is the consequence of hyperbolicity of wave equations
and  locality of the interaction, see, for instance,
\cite[ch. 5, \S 28]{Vla66},
Reed, Simon
\cite[v.2, Theorem X.77]{ReedSimon75}.

\medskip
To construct the required time column, and/or
the corresponding
light cone, we formulate as Theorem 7.2
 the assertion proved by
S.~Nelson
\cite{SNelson71}.

\medskip
{\large Theorem 7.2}
 (see
\cite[Corollary 2]{SNelson71},
our notations follow the notation in
\cite{SNelson71},
 see also the required conditions on
the considered solutions
\cite{SNelson71},
\cite{Pan80-82}).
  {\it Let  $u_{in}$ be a free solution that
satisfies the conditions
$(6.8).$
Then
$$
\lim_{t\to+\infty} t^{3/2} e^{i\alpha(t,\lambda)}
u^+_{in}(t,{\lambda t\over\sqrt{1+\lambda^2}})
= (1+\lambda^2)^{5/4}
(u^+_{in})^\sim(\lambda),
$$
where $\lambda \in \R^3,$
 $$\alpha(t,\lambda) \equiv {3\pi\over 4} +
{t\over \sqrt{1 + \lambda^2}}$$ and $u^+_{in},$
$(u^+_{in})^\sim$
is the  positive frequency part of the free
(real) solution at the time zero and, correspondingly,
its Fourier transform.
}

\medskip
{\it Remark.}
   S.~Nelson
\cite{SNelson71}
uses the choice of space variables,  corresponding
to the choice of space variables in the form
 $(mt,mx).$ In momentum space this choice corresponds
to the choice of coordinates in the form  $(k^0/m,{\bf
k}/m).$  Here $m$ is the mass constant from the nonlinear
equation.

\medskip
We continue the proof of Theorem  7.1.

We apply now Theorem  7.2 (see
\cite[Corollary 2]{SNelson71})
to the initial $in$-\-data.
For this purpose we construct a column in which
 $u(t,x) \not = 0.$ This column is constructed
with the help of a point at which
$(u^+_{in})^\sim(b) \not = 0$
and a Lorentz rotation.

We note that we need the assertion about vanishing
in the form of
  Theorem 3.1
 about unitarizability.
The conditions of Theorem 3.1 imply that
$(u^+_{in})^\sim(k)$ is continuous in $k.$
Since  $u^+_{in} \not
= 0,$
so there is a point
 $b$
in momentum space
such that
$(u^+_{in})^\sim(b) \not = 0.$
This fact is equivalent to the existence
of a Lorentz transformation such that
$(u^+_{in})^\sim(0) \not = 0.$
Here
$
u_{in,\Lambda}(t,x) = u_{in}(\Lambda(t,x))
$
and
$$
(u^+_{in,\Lambda})^\sim(k)=(2\pi)^{-3/2}
\int e^{ikx} u_{in,\Lambda}
(0,x) d^3 x.
$$
Let $c(u_{in}) := |(u^+_{in,\Lambda})^\sim(0)| > 0.$
There exists such  $r_1,$ that
$$
\sup_{|k| \leq r_1} |(u^+_{in,\Lambda})^\sim(k)
-
(u^+_{in,\Lambda})^\sim(0)|
\leq {1\over 8} c(u_{in}),
$$
$$
\sup_{|k|\leq r_1}((1+{k^2\over m^2})^{5/4}-1)
\leq {1\over 8} c(u_{in}).
$$
Let
$$
u^+_{as,\Lambda}(t,x) = |mt|^{-3/2}(1 +
{x^2\over \sqrt{t^2-x^2}})^{5/4}
e^{i\alpha(t,x)}
(u^+_{in,\Lambda})^\sim({x\over\sqrt{t^2-x^2}}),
$$
here
 $$\alpha(t,x) = {3\pi\over 4} + m\sqrt{t^2-x^2}$$
and we have used the transform
$x = \lambda t/\sqrt{m^2 + \lambda^2},$
that is,
$\lambda = m x/\sqrt{t^2-x^2},$
to introduce some approximation.
This  approximation allows to use the assertion of
Theorem  7.2.
We note that
$\mbox{Re }u^+_{as,\Lambda}(t,x)$ corresponds to
the approximation of the initial solution $u_{in}(t,x).$
We write
\be
u_\Lambda(t,x)
&=&
\mbox{ Re }u^+_{as,\Lambda}(t,x)
 + u_{in,\Lambda}(t,x) -\mbox{ Re }u^+_{as,\Lambda}(t,x)
\\
&&+ u_\Lambda(t,x)-u_{in,\Lambda}(t,x).
\ee
We take $r$ (this $r$ defines the width of the column) and
we suppose that
$|x| < r.$ Moreover, this value of $r$ we choose
sufficiently  large with respect to chinks (of the size
$\leq 2\pi/\it{mass\;\ constant}$) between the constructed
cylinders and such that the causal envelope of
constructed cylinders contains the column with the
width $r/2.$

There exists such (sufficiently large)
$a_1 > 0,$
that for all $t < -a_1,$ $|x| < r$
$$
\sup_{|x|\leq r} |mt|^{3/2}\, \vert u_{in,\Lambda}(t,x)
- \mbox{Re }u^+_{as,\Lambda}(t,x) \vert
\leq{1\over 8}c(u_{in}),
\eqno(7.4)
$$
$$
|mt|^{3/2} \| u_\Lambda(t,x)
- u_{in,\Lambda}(t,x)\|_\infty
\leq {1\over 8}c(u_{in}),\eqno(7.5)
$$
(7.4) is implied by  Theorem 7.2, and (7.5) follows from the
requirements on norms of initial   $in$-data
\footnote{We note that $t^{3/2}$ convergence is the
consequence of conditions that the solution $u_{in}$ satisfied
or the consequence of the same conditions for the Lorentz
rotated solution
$u_{in,\Lambda}.$
Here the Lorentz rotation is defined by
(non-zero) solution
 $u_{in}$ and the requirement
$(u^+_{in,\Lambda})^\sim(0) \not = 0$ (under conditions
that are satisfied by the solution $u_{in}$). In any case these
conditions are fulfilled with our choice of initial data from
$\Sc_{Re}(\R^3)\oplus\Sc_{Re}(\R^3)$ due to
the choice of sufficiently small $\alpha$ in Theorem 3.1.

The choice of the Lorentz rotation (or the choice of the
Lorentz frame) and its consideration is equivalent to the
fact that instead of the solution $u(t,x)$
we consider the solution
$u_\Lambda(t,x) = u(\Lambda(t,x))$
with the straight forwarded column, directed into the past.}.
  For $t_n\leq t\leq t_n+\pi/2m,$ where
$$
 t_n
 = {1\over m}({\pi \over 2} +
 \pi n) - \arg(u_{in}^+)^\sim(0)
$$
and $m$ is the mass constant, the explicit form of
$u^+_{as,\Lambda}(t,x)$
implies that
$$
\vert\mbox{\,Re }u_{as,\Lambda}(t,x) \vert \geq
|mt|^{-3/2} {1\over 8}c(u_{in})>0.
$$
Here we require that  $a_1$ is larger than   $4m a^2/\pi.$
Then, we choose  $r,$ defined the width of the cylinder,
such that it is larger than the chinks between the cylinders
and such that the causal envelope of these cylinders contains
some (more narrow) column (of nonzero width).

 The causal envelope of the
 cylinder contains the part of some (light) cone
from the top to the section of the cone.
 The section of this cone coincides
  with the support of the
cylinder.
 This is implied  by hyperbolicity,
 the fact that $u(t,x)$ for $t\leq 0$
satisfies the free wave equation and, therefore, the solution
 $v_{\Lambda}(t,x),$ corresponding to  the consideration of
 $u_{\Lambda}(t,x)$
 (i.e. corresponding to the consideration
of the solution in some Lorentz frame), satisfies the free
wave equation in some  domain. This domain is  equal to the
domain
$\{(t,x)\in\R^4 \,\vert \,t\leq 0\}$ in this Lorentz frame.
 This domain with the space-like boundary contains the light
cone with the top at zero (in the considered Lorentz frame)
and directed into the past. Therefore, the causal envelope
of (the set of) cylinders contains the column and, thus,
it is the light cone directed into the past.

 Therefore, causality and hyperbolicity  imply that
$v(t,x)=0$ in the backward light cone with the top
in the sufficiently large past. Consequently,
 $v=0,$ see
\cite[Corollary 1, Corollary $1^\prime$]{SNelson71},
\cite{Mor63}.
We note that by S.~Nelson
\cite[Corollary 1, Corollary $1^\prime$]{SNelson71}
  Corollary 1, as far as Corollary $1^\prime$, implies that
$$
\Vert v^+\Vert_{L_2} =
\lim_t \Vert v(t,\cdot)\Vert_{L_2}
=
\lim\Vert\chi_\Gamma v(t,\cdot)\Vert_{L_2}.
$$
Here $\chi_\Gamma$ is the characteristic function of the
cone $\Gamma.$
Since we consider the solutions $v$ with the
initial data in
 $H^{1/2}(\R^3)\oplus H^{-1/2}(\R^3),$ it is more
simple to use the S.~Nelson assertion
\cite{SNelson71},
or even the corresponding assertion in
the Vladimirov book
\cite[ch. V, \S 29, Theorem 1, Corollary 2]{Vla66}.
{\it Theorem $7.1$ is proved.}



\medskip
\medskip
\medskip
\medskip
{\large Appendix}
\medskip
\medskip

As pointed out by P. Kumlin
\footnote{We are indebted to P. Kumlin for sending us this
improvement.}
 Lemma 3.4 in the
proof of real analyticity
\cite{Kumlin92}
is not correct as it is
formulated.  That is also the case for Proposition 1.4 in
Brenner
\cite{Brenner88-89},
which is the origin of Lemma 3.4.  However the conclusion in
Step 3 remains true if lines 4 to 25 on page 265
\cite{Kumlin92}
are replaced by
the argument below.

It remains to show that, for all
$t \in [-T,T],$
$$
{\cal K}(t) := \{\int_{0}^{t}
K(t-s)(\varphi^2\chi_R\psi_j(s))ds : j = 1,2,3,...\}
$$
has a convergent subsequence in $L_6.$
 By the Rellich-\-Kon\-dra\-shov
theorem it suffices to prove that
$
{\cal K}(t)$
 is bounded in some $L_{p'}^{s'+\epsilon}, \epsilon >0, $
where $
{1\over {p'}}-{{s'}\over 3} = {1\over 6}
$ and $
s'>0.$
Set $p={3\over 2}$
and $
s=
{1\over 6}+{\epsilon\over 3}$
 $(s' = {1\over 2}).$
For a fixed $t \in [-T,T],$ Proposition 3.1 yields
$$
\| \int_0^t K(t-s)
(\varphi^2\chi_R\psi_j(s))ds \|_{
L_{p'}^{s'+\epsilon}} \le
 \int_0^t
k(t) \|
\varphi^2\chi_R\psi_j(s)
\|_{L_{p}^{s}} ds$$
 $$ \le C \sup_{s\in [0,t] }
\|\varphi^2\chi_R\psi_j(s)
\|_{L_{p}^{s}} .$$
Here we introduce the Besov spaces
 $
B = B
_{p}^{s+\epsilon, q} $ with norm
$
(0<s<1) $
$$
\|f\|_B=(
\int_0^\infty
(t^{-s}
\omega_p(f;t))^q{dt\over t})^{1/q}.
$$
$\omega_p(f;t)$ denote the continuity modulus
$$\omega_p(f;t)
=\sup_{|\nu|\le t}\|f(\cdot +\nu)-f(\cdot)\|_{L_p}$$
with the usual modification for
$
p=\infty.$
The embedding
$L_{p}^{s+2\epsilon} \subset
B_{p}^{s+\epsilon,q} \subset
L_{p}^{s},
$ $1\le q < \infty$ is well known.
Straight forward calculations give
$$
\omega_p
(\varphi^2\chi_R\psi_j;t)
\le \|\chi_R\|_{L_\infty}
\{\|\varphi\|_{L_6}^2
\omega_3(\psi_j;t)+...\}
$$
$$ +C\|\varphi\|_{L_6}
\omega_\infty(\chi_R;t),
$$ where $...$ denotes cyclic permutation of
$\varphi,\varphi $ and $       \psi_j.$
 By embedding above and Sobolev's inequality it follows that
$$\sup_{x\in {\bf R} }
\|\varphi^2\chi_R\psi_j(s)
\|_{B_{p}^{s+\epsilon,1}}
\le C
\|\varphi\|_{{\bf Z}}^2\|\psi_j\|_{{\bf Z}}
$$
for $
\epsilon $
 small enough and we conclude
 $$
\| \int_0^t K(t-s)
(\varphi^2\chi_R\psi_j(s))ds \|_{
L_{p'}^{s'+\epsilon}} \le
C \sup_{s\in[0,t]}\|\varphi^2\chi_R\psi_j(s)
\|_{B_{p}^{s+\epsilon,1}}$$
 $$ \le C \|
\varphi\|_{{\bf Z}}^2\|\psi_j\|_{{\bf Z}}
$$


\medskip
\medskip
\medskip
\medskip
{\large Acknowledgments}
\medskip
\medskip

This is the first paper  of the project
 $\phi^4_4 \cap M.$
The one of the goal of this project is to support
partly the Russian Fundamental Research Foundation.

We acknowledge
Anatoly Kopilov, Valery Serbo,
 Vasily Serebrjakov,
 Ludwig Faddeev, Anatoly Vershik,
 Peter Osipov,  Volja Heifets, Zinaida and Julia
 for the help, advice, and criticism.


\end{document}